\documentclass[11pt,a4paper]{article}
\usepackage[T1]{fontenc}
\usepackage[latin1]{inputenc}
\usepackage[francais,english]{babel}
\usepackage[top=4cm, bottom=4cm, left=2.5cm, right=2.5cm]{geometry}
\usepackage{amsmath,amsthm,mathrsfs,amssymb,amsfonts,nicefrac,stmaryrd,yhmath}
\usepackage{dsfont}
\usepackage{shadethm}
\usepackage{graphicx}
\usepackage[dvips,svgnames]{xcolor}
\usepackage{enumerate}
\usepackage{array}
\usepackage[normal]{multicap}
\usepackage{pstricks,pstricks-add,pst-math,pst-xkey}
\usepackage{pstricks}
\usepackage{pst-node,pst-3d}
\usepackage{varioref}
\usepackage{fancyhdr,fancybox,shadow}
\usepackage[normalem]{ulem}
\usepackage{color}
\renewcommand{\textbf}[1]{{\bfseries\mathversion{bold}#1}}
\newcommand{\R}{\mathbb{R}}
\newcommand{\C}{\mathbb{C}}

\newcommand{\N}{\mathbb{N}}

\newenvironment{disarray}{\everymath{\displaystyle\everymath{}}\array}{\endarray}
\newcommand{\dis}{\displaystyle}
\newcommand{\abs}[1]{\left|#1\right|}
\newcommand{\eps}{\varepsilon}
\newcommand{\norme}[1]{\left\|#1\right\|}
\renewcommand{\a}{\alpha}

\renewcommand{\leq}{\leqslant}
\renewcommand{\geq}{\geqslant}

\renewcommand{\tild}{\widetilde}
\newcommand{\re}{\Re\mathrm{e}}
\newcommand{\im}{\Im\mathrm{m}}

\theoremstyle{plain}

\newtheorem{theo}{Theorem}[section]

\newtheorem{prop}{Proposition}[section]

\newtheorem{lemme}{Lemma}[section]
\newtheorem{cor}{Corollary}[section]

\newenvironment{dem}{\noindent\textit{Proof :}}{\hfill\rule{2mm}{2mm}}

\numberwithin{equation}{section}

\begin{document}
\title{Transition between linear and exponential propagation in Fisher-KPP type reaction-diffusion equations}
\author{Anne-Charline COULON$^{\small{*}}$, Jean-Michel ROQUEJOFFRE$^{\small{*}}$ \\
\footnotesize{$^{\small{*}}$Institut de Mathématiques, (UMR CNRS 5219), Université Paul Sabatier,}\\
\footnotesize{118 Route de Narbonne, 31062 Toulouse Cedex, France}}
\date{}
\maketitle
\begin{abstract}
We study the Fisher-KPP equation with a fractional laplacian of order $\a \in (0,1)$. We know that the stable state invades the unstable one at constant speed for $\a = 1$, and at an exponential in time velocity for $\a \in (0,1)$. The transition between these two different speeds is examined in this paper. We prove that during a time of the order $-\ln(1-\a)$, the propagation is linear and then it is exponential.
\end{abstract}
\section{Introduction}

We are interested in the large time behaviour of the solution $u$ to the evolution problem:
\begin{equation} \label{sys}
u_t +(-\Delta)^{\a}u=u-u^2, \quad \R^d, t>0.
\end{equation}
The nonlinearity $u-u^2$ is often referred to a Fisher-KPP nonlinearity, from the reference
 \cite{Kolmo}, and is motivated by spatial propagation or spreading of biological species. Traditionally, two kinds of data are considered: first we study a compactly supported initial data. It corresponds to the study of spatial spreading when at initial time some areas are not invaded by the population. Then we consider a nondecreasing initial data, in one space dimension. This choice is motivated by \cite{Kolmo}.
In this paper, we study the transition between two speeds of propagation:
\begin{itemize}
\item When $\a$ is equal to $1$, it is well known (see \cite{AW}) that the stable state invades the unstable one at constant speed equal to 2: there is a linear propagation. Many papers deal with this phenomenon. Let us mention one of the most general works  \cite{BHN}: here the authors introduce a very general notion of propagation velocity and prove it can behave in various ways in general unbounded domains of $\R^d$. For instance, for some locally $\mathcal{C}^2$ domains of $\R^d$ which do not satisfy the extension property (that is to say the boundary can be covered by a sequence of open sets whose radii do not tend to $0$), the spreading speed may be infinite, and for some locally $\mathcal{C}^2$ domains of $\R^d$ which satisfy the extension property and are strongly unbounded in every direction of $\mathbb{S}^1$, the spreading speed may be zero. See also \cite{BHNI} and \cite{Wein} for spatially periodic media.

\item When $\a \in (0,1)$, the authors of \cite{BRR} prove there is still invasion of the unstable state by the stable state, and in \cite{pre} it is proved that propagation holds at an exponential in time velocity.  This result provide a mathematically rigorous justification of numerous heuristics about this model. (see \cite{vulpi} for instance) Note that the position of the level sets of $u$ move exponentially fast as $t \mapsto + \infty$ for integro-differential equations for example. (see \cite{Jimmy}). Also note that exponentially propagating solutions exist in the standard KPP equations, as soon as the initial datum decays algebraically. This fact was noticed by Hamel and Roques in \cite{HR}.
\end{itemize}
So, we want to understand, in a more precise fashion, the transition between the two models.

First, we study the fundamental solution $p$ to \eqref{sys}, for $d\geq 1$. Indeed, it has a special role and all the heuristics are based on it. As a matter of fact, all the formal asymptotic studies (see for instance \cite{vulpi} or \cite{Cas}) are based on the study of the fundamental solution. We will see that it has the following expansion as $t>0$ and $\a$ tends to $ 1$:
$$
\begin{disarray}{rclr}
\abs{p(x,t)-\frac{C_{\a}\sin(\a \pi)t}{\abs{x}^{d+2\a}} -\frac{e^{-\frac{\abs{x}^{2 \a}}{4t}} }{(4 \pi t)^{\nicefrac{d}{2}}  \abs{x}^{d(1-\a)}}} & \leq & C  \frac{(1-\a)t^2} { \abs{x}^{d+4\a}}, & \forall x \in \R^*,
\end{disarray}
$$
where $C_{\a}$ and $C$ are positive constants.
The study of this inequality reveals a critical time scale $\tau_{\a} $ of the order $ -  \ln(1-\a)$ for $\a$ close to $1$, where the transition occurs. 
This leads to the main theorems of this paper, which parallel the main results of the preprint \cite{pre}:
\begin{theo} \label{thm}
Consider $u$ the solution to \eqref{sys}. Let $u_0$ be a function with compact support, $0\leq u_0 \leq 1$, $u_0 \neq 0$. Set $\tau_{\a} =  -\ln (1-\a)$. Then:
\begin{itemize}
\item For all $\sigma > 2$, $u(x,t) \rightarrow 0$ uniformly in $\{ \abs x \geq \sigma t^{\nicefrac{1}{\a}}\}$ as $\a \rightarrow 1, t \rightarrow +\infty, t < \tau_{\a}$.
\item For all  $0< \sigma < 2$, $u(x,t) \rightarrow 1$ uniformly in $\{ \abs x \leq \sigma t^{\nicefrac{1}{\a}}\}$ as $\a \rightarrow 1, t \rightarrow +\infty, t<\tau_{\a}$,
\end{itemize}
Moreover, there exists a constant $C>0$ independant of $\a$ such that:
\begin{itemize}
\item For all  $\sigma > \dis{\frac{1}{d+2\a}}$, $u(x,t) \rightarrow 0$ uniformly in $\{ \abs x \geq e^{\sigma t}\}$ as $\a \rightarrow 1, t \rightarrow +\infty, t > C \tau_{\a}$.
\item For all  $0< \sigma < \dis{\frac{1}{d+2\a}}$, $u(x,t) \rightarrow 1$ uniformly in $\{ \abs x \leq e^{\sigma t}\}$ as $\a \rightarrow 1, t \rightarrow +\infty, t >C  \tau_{\a}$.
\end{itemize}
\end{theo}

The second theorem concerns initial data that are increasing in space:
\begin{theo} \label{thm2}
Consider $u$ the solution to \eqref{sys}. Let $u_0$ be a function  in $[0,1]$  measurable, nondecreasing, such that $\underset{x \mapsto +\infty}{\lim}u_0(x)=1$ and satisfying $$u_0(x) \leq c e^{-\abs x ^{\a}}, \mbox{for x $\in \R_-$ },  $$ for some constant $c$. Set $\tau_{\a} =  -\ln (1-\a)$. Then:
\begin{itemize}
\item For all $\sigma > 2$, $u(x,t) \rightarrow 0$ uniformly in $\{  x \leq -\sigma t^{\nicefrac{1}{\a}}\}$ as $\a \rightarrow 1, t \rightarrow +\infty, t < \tau_{\a}$.
\item For all  $0< \sigma < 2$, $u(x,t) \rightarrow 1$ uniformly in $\{  x \geq -\sigma t^{\nicefrac{1}{\a}}\}$ as $\a \rightarrow 1, t \rightarrow +\infty, t<\tau_{\a}$,
\end{itemize}
and there exists a constant $\overline{C} >0$ independent of $\a$ such that:
\begin{itemize}
\item if $\sigma > \dis{\frac{1}{2\a}}$, $u(x,t) \rightarrow 0$ uniformly in $\{  x \leq -e^{\sigma t}\}$ as $\a \rightarrow 1, t \rightarrow +\infty, t > \overline{C} \tau_{\a}$.
\item if $0< \sigma < \dis{\frac{1}{2\a}}$, $u(x,t) \rightarrow 1$ uniformly in $\{ x \geq -e^{\sigma t}\}$ as $\a \rightarrow 1, t \rightarrow +\infty, t >\overline{C}  \tau_{\a}$.
\end{itemize}
\end{theo}
For $t \in [\tau_{\a},C \tau_{\a}]$, there is a transition between two different speeds: this phenomenon will not be studied in this paper.\\
\textbf{Remark 1:}  For $1 \leq t\leq \tau_{\a}$, we have: $ t- t^{\nicefrac{1}{\a}} = \underset{\a \rightarrow 1}{\cal O} (\sqrt{ 1-\a }) $. So, we have $\sigma t^{\nicefrac{1}{\a}} \underset{\a \rightarrow 1}{\sim} \sigma t$ in the range $t \in [0, \tau_{\a}]$ and so the propagation is truly linear.\\
\textbf{Remark 2:} The decay $e^{-\abs x ^{\a}}$ is almost optimal. Indeed, when $\a=1$, recall that \cite{AW} implies linear propagation at velocity $2$. Now, if $u_0(x) \underset{x \mapsto - \infty}{\sim} e^{-\gamma\abs x }$, with $\gamma <1$, \cite{uchi} implies linear propagation at velocity $\frac{1}{\gamma}+\gamma$. This illustrates the fact that, in order to obtain a result that is uniform in $\a$, we need the $e^{-\abs x ^{\a}}$ decay at least.\\
\textbf{Remark 3:} Adapting the proof of theorem \ref{thm}, we could prove that with an initial datum of the form $e^{-\gamma \abs x ^{\a}}$, $\gamma <1$, there is a linear propagation phase with velocity $\frac{1}{\gamma}+\gamma$ in a time interval of length $\tau_{\a}$.

The proofs of the main theorems \ref{thm} and \ref{thm2} follow the same plan. An upper bound is obtained thanks to an appropriate estimate of the fundamental solution. Concerning the lower bound, the study of the fundamental solution is not enough, we need an iterative scheme, and in each iteration the solution has to be suitably truncated. It enables us to know the evolution of the level set of the form  $\{ x \in \R^d \  | \ u(x,t)=(1-\a)^{1+\kappa}, \kappa >0 \}$. However, we need to study level sets that do not depend on $\a$: an intermediate proposition will ensure the connection between these two level sets.

The paper is organized as follows. Section 2 contains the study of the fundamental solution. In Section 3 we prove the intermediate result that makes possible the connection between levels sets that depend on $\a$ and the one that do not depend on $\a$. Section 4 and 5 are concerned with the proof of the main theorems of the paper, concerning respectively compactly supported initial data and nonincreasing initial data.

\textit{Notation:} Throughout this article $C$ denotes, as usual, a constant independent of $\a$.

\textit{Acknowledgement:} Both authors were supported by the ANR project PREFERED. They thank Professor X. Cabr\'e for fruitful discussions.

\section{The fundamental solution}\label{2}
First, for the sake of simplicity, we consider the one space dimension case to underline the idea of the proof. The higher space dimension case is treated in subsection \ref{D} and requires the use of special functions like the Bessel and Whittaker functions.

\subsection{In one space dimension}

Recall that we are solving:
\begin{equation} 
\left\{
\begin{array}{rcl}
p_t +(-\partial_{xx})^{\a}p&=&0, \quad \R, t>0\\
p(0,x)&=&\delta_0(x), \hfill x \in \R, \\
\end{array}
\right.
\end{equation}
whose solution is: $p(x,t)=\mathcal{F}^{-1}(e^{-\abs{\xi}^{2\a}t})= t^{-\nicefrac{1}{2\a}} p_{\a}(x t^{-\nicefrac{1}{2\a}})$, where:

$$p_{\alpha}(x)= \frac{1}{2 \pi}\int_{\R} e^{-\abs \xi^{2 \a}}e^{i x. \xi} d\xi.$$

\begin{prop}\label{in1D}
We have: $$
\begin{disarray}{lrclr}
&\abs{p_{\alpha}(x)-\frac{\Gamma(2\a+1) \sin(\a \pi)}{ \pi \abs{x}^{1+2\a}} -\frac{e^{-\frac{\abs{x}^{2 \a}}{4}} }{2 \sqrt{\pi}  \abs{x}^{1-\a}} }& \leq & C  \frac{(1-\a)} { \pi \abs{x}^{1+4\a}}, & \forall x \in \R^* \\
\mbox{ and } & \ \ \abs{p_{\a}(x)-p_{\a}(0)} &\leq& \frac{\Gamma\left(\nicefrac{1}{\a}+1\right)}{2\pi} \abs x , & \forall x \in \R  \\
\end{disarray}
$$
\end{prop}
The first inequality will be used to control $p_{\a}$ for large values, whereas the second one will be used to control $p_{\a}$ in the vicinity of $0$. The proof is based on \cite{Polya} and \cite{Blu}.

\begin{dem}
By symmetry we have, for all $x \in \R$:
$$
\begin{disarray}{rcl}
p_{\alpha}(x)
&=& \frac{1}{ \pi}\int_{0}^{\infty} e^{-r^{2 \a}} \cos(r \abs x) dr\\
\end{disarray}
$$
To prove the first inequality, we take $x \ne 0$.
Then, we integrate by parts and introduce the new variable $u=\abs{x}^{2\a} r^{2 \a}$  as in \cite{Polya} and \cite{Blu}. Thus:
\begin{equation} 
p_{\alpha}(x)
= \frac{1}{ \pi \abs{x}^{1+2\a}}\int_{0}^{\infty} e^{- u \abs{x}^{-2 \a}} \sin(u^{\frac{1}{2\a}}) du \\
\end{equation}
Now, we have to study: $I_{\a}(x) := \dis{\int_{0}^{\infty} } e^{- u \abs{x}^{-2 \a}} e^{i u^{\frac{1}{2\a}}} du$. We rotate the line of integration by $\frac{\pi}{2}$:
$$
\begin{disarray}{rcl}
I_{\a}(x)&=& \int_0^{\infty}i e^{-r^{\nicefrac{1}{2\a}}\sin(\frac{\pi}{4\a})} e^{i r^{\nicefrac{1}{2\a}}\cos(\frac{\pi}{4\a})} e^{-ri\abs{x}^{-2\a}} dr
\end{disarray}
$$
\noindent First we denote by: 
$$I_{\a, \infty}= \int_0^{\infty}i e^{-r^{\nicefrac{1}{2\a}}\sin(\frac{\pi}{4\a})} e^{i r^{\nicefrac{1}{2\a}}\cos(\frac{\pi}{4\a})}  dr,$$
  which can be simplified by rotating the interval of integration  by $\a \pi - \frac{\pi}{2}$:
$$
\begin{disarray}{rclcl}
I_{\a, \infty}&=&\int_0^{\infty}e^{-u^{\frac{1}{2\a}}}e^{i \a \pi} du&=& \Gamma (2\a+1) e^{i \a \pi}\\
\end{disarray}
$$
Secondly, we write:
$$I_{\a}^1(x):= \int_0^{\infty}i e^{-r^{\nicefrac{1}{2\a}}\sin(\frac{\pi}{4\a})} e^{i r^{\nicefrac{1}{2\a}}\cos(\frac{\pi}{4\a})}( e^{-ri\abs{x}^{-2\a}} -1) dr,$$
so that: $I_{\a}(x)=I_{\a, \infty}+I_{\a}^1(x)$.
Then we introduce: $f_r(\a)=e^{-r^{\frac{1}{2\a}}\sin(\frac{\pi}{4\a})} e^{i r^{\frac{1}{2\a}}\cos(\frac{\pi}{4\a})}$, and write:
$$I_{\a}^1(x)=  \int_0^{\infty}i f_r(1) e^{-ri\abs{x}^{-2\a}} dr - \int_0^{\infty}i f_r(1) dr+\int_0^{\infty}i ( f_r(1)-f_r(\a)) ( e^{-ri\abs{x}^{-2\a}} -1) dr.$$
Taking the imaginary part of each element of the right hand side:
\begin{itemize}
\item  $\begin{disarray}{rcl} \int_0^{\infty}i f_r(1) dr
&=&  \int_0^{\infty}e^{-u^{\nicefrac{1}{2}}}du,
\end{disarray}$
so: $\im \left(\displaystyle{\int_0^{\infty}} i f_r(1) dr \right)=0$
\item $ \displaystyle{ \int_0^{\infty}} i f_r(1) e^{-ri\abs{x}^{-2\a}}  dr$ is treated by rotating back the  integration line by $-\frac{\pi}{2}$. Then:
$$
\begin{disarray}{rcl}
\im \left(\int_0^{\infty}i f_r(1) e^{-ri\abs{x}^{-2\a}}  dr \right)&=&\im \left(\int_0^{\infty}e^{iu^{\frac{1}{2}}-u\abs{x}^{-2\a}} du\right)\\
&=& \pi \abs{x}^{3\a} p_{1}(x^{\a}) \\
&=&\frac{\sqrt{\pi}}{2 } e^{-\frac{\abs{x}}{4}^{2 \a}} \abs{x}^{3 \a} \\
\end{disarray}
$$
\item $\displaystyle{\int_0^{\infty}} i ( f_r(1)-f_r(\a)) ( e^{-ri\abs{x}^{-2\a}} -1) dr$ is of lower order since it is bounded by $C(1-\a) \abs{x}^{-2\a}$, where $C$  is a constant independent of $\a$.
Indeed: $$
\begin{disarray}{rcl}
\abs{\im \left(\int_0^{\infty}i ( f_r(1)-f_r(\a)) ( e^{-ri\abs{x}^{-2\a}} -1) dr \right)}& \leq & \int_0^{\infty} \abs{f_r(1)-f_r(\a)} \abs{e^{-ri\abs{x}^{-2\a}} -1} dr \\
& \leq & \int_0^{\infty}(1- \a) \sup_{y\in (\a,1)} \abs{\partial_y f_r(y)} r \abs{x}^{-2\a}  dr
\end{disarray}
$$
Moreover, for $y \in (\a,1)$:
$$
\begin{disarray}{rcl}
\abs{\partial_y f_r(y)}&=&\abs{\left( -i \frac{ \ln(r)}{2 y^2}+\frac{\pi}{4y^2}\right)  r^{\frac{1}{2y}} e^{i\frac{\pi}{4y}}e^{-r^{\frac{1}{2y}}\sin(\frac{\pi}{4y})} e^{i r^{\frac{1}{2y}}\cos(\frac{\pi}{4y})}}\\
& \leq & \left( \frac{\pi}{4y^2}+  \frac{\abs{ \ln(r)}}{2 y^2} \right) r^{\frac{1}{2y}} e^{-r^{\frac{1}{2y}}\sin(\frac{\pi}{4y})}\\
& \leq &\left(\frac{\pi}{2} - \ln (r)\right) \frac{r^{\frac{1}{2}}}{2 \a^2} e^{-\frac{r^{\frac{1}{2\a}}}{\sqrt{2}}} \mathds{1}_{(0,1)}(r) + \left(\frac{\pi}{2} + \ln (r)\right)\frac{r^{\frac{1}{2\a}}}{2 \a^2} e^{-\frac{r^{\frac{1}{2}}}{\sqrt{2}}} \mathds{1}_{(1,+\infty)}(r)\\
& \leq &  \left(\frac{\pi}{2} - \ln (r) \right)\frac{1}{2 \a^2}  \mathds{1}_{(0,1)}(r) +\left(\frac{\pi}{2} + \ln (r)\right)\frac{r^{\frac{1}{2\a}}}{2 \a^2} e^{-\frac{r^{\frac{1}{2}}}{\sqrt{2}}} \mathds{1}_{(1,+\infty)}(r)\\
\end{disarray}
$$
Thus, for $\a \geq \frac{1}{2}$:
$$
\begin{disarray}{rcl}
\multicolumn{3}{l}{\abs{\im \left(\int_0^{\infty}i ( f_r(1)-f_r(\a)) ( e^{-ri\abs{x}^{-2\a}} -1) dr \right)}}\\
 &\leq & \int_0^1 (1-\a) \abs{x}^{-2 \a} r  \left(\frac{\pi}{2} - \ln (r)\right)\frac{1}{2 \a^2} dr + \int_1^{\infty} (1-\a) \abs{x}^{-2 \a} \left(\frac{\pi}{2} + \ln (r)\right)\frac{r^{\frac{1}{2\a}+1}}{2 \a^2} e^{-\frac{r^{\frac{1}{2}}}{\sqrt{2}}}dr \\
& \leq &C (1-\a) \abs{x}^{-2 \a}+C (1-\a) \abs{x}^{-2 \a}\left( \int_1^{\infty} r^2 e^{-\frac{r^{\frac{1}{2}}}{\sqrt{2}}}dr+\int_1^{\infty}\ln(r)  r^2 e^{-\frac{r^{\frac{1}{2}}}{\sqrt{2}}}dr\right)\\
& \leq & C  (1-\a) \abs{x}^{-2 \a}.\\
\end{disarray}
$$
\end{itemize}
Consequently:
$$
\begin{disarray}{rcl}
p_{\alpha}(x)&=&\frac{1}{ \pi \abs{x}^{1+2\a}} \im\big(\int_{0}^{\infty} e^{- u \abs{x}^{-2 \a}} e^{i u^{\frac{1}{2\a}}} du\big)\\
&=&\frac{\Gamma(2\a+1) \sin(\a \pi)}{ \pi \abs{x}^{1+2\a}} +\frac{e^{-\frac{\abs{x}}{4}^{2 \a}} \abs{x}^{3 \a}}{2 \sqrt{\pi} \abs{x}^{1+2\a}} + \frac{ \im \left(\int_0^{\infty}i ( f_u(1)-f_u(\a)) ( e^{-r i \abs{x}^{-2\a}} -1) dr \right)}{ \pi \abs{x}^{1+2\a}}\\
\end{disarray}
$$
Thus we obtain:
\begin{eqnarray}
\abs{p_{\alpha}(x)-\frac{\Gamma(2\a+1) \sin(\a \pi)}{ \pi \abs{x}^{1+2\a}} -\frac{e^{-\frac{\abs{x}}{4}^{2 \a}} }{2 \sqrt{\pi}  \abs{x}^{1-\a}} }& \leq & C  \frac{(1-\a)} { \pi \abs{x}^{1+4\a}}, \  \forall x \in \R^*.
\end{eqnarray}

To get an estimate for $x$ in a compact set, we use the expression of $p_{\a}(0)$ given by:
$$
\begin{disarray}{rcl}
p_{\a}(0) =  \frac{1}{2 \pi}\int_{\R} e^{-\abs \xi^{2 \a}} d\xi = \frac{\Gamma(\nicefrac{1}{(2 \a)}+1) }{\pi} < {+\infty}
\end{disarray}
$$
For all $y \in \R$, we have:
$$
\begin{disarray}{rclcl}
 \abs{p_{\a}'(y)} & \leq& \frac{1}{2\pi} \int_{\R} \abs {\xi} e^{-\abs{\xi}^{2\a}}d\xi &=&\frac{1}{2\pi} \Gamma\left(\nicefrac{1}{\a}+1\right)\\
\end{disarray}
$$
So, by the mean value theorem:
$$
\begin{disarray}{rclr}
 \abs{p_{\a}(x)-p_{\a}(0)} &\leq& \frac{\Gamma\left(\frac{1}{\a}+1\right)}{2\pi}\abs x  , & \forall x \in \R  \\
\end{disarray}
$$
\end{dem}
\subsection{In higher space dimension}\label{D}
In this case, we are solving:
\begin{equation*} 
\left\{
\begin{array}{rcl}
p_t +(-\Delta)^{\a}p&=&0, \quad \R^d, t>0\\
p(0,x)&=&\delta_0(x), \hfill x \in \R^d, \\
\end{array}
\right.
\end{equation*}
whose solution is: $p(x,t)=\mathcal{F}^{-1}(e^{-\abs{\xi}^{2\a}t})= t^{-\nicefrac{d}{2\a}} p_{\a}(x t^{-\nicefrac{1}{2\a}})$, where:
$$p_{\alpha}(x)= \frac{1}{(2 \pi)^d}\int_{\R^d} e^{-\abs \xi^{2 \a}}e^{i x. \xi} d\xi.$$
\begin{prop}\label{inD}
We have: 
$$
\begin{disarray}{lrclr}
&\abs{p_{\a}(x) - \frac{2(2\pi)^{-\frac{d+1}{2}}\sin(\a \pi) D_{\a}}{\abs x^{d+2\a}}- \frac{ e^{-\frac{\abs x }{4}^{2\a}}}{(4 \pi)^{\frac{d}{2}} \abs x^{(1-\a)d}\a}}&\leq& \frac{C(1-\a)}{\abs x^{d+4\a}}, & \forall  x \in \R^*\\
\mbox{ and } &\abs{p_{\a}(x)-p_{\a}(0)} &\leq&  \frac{\Gamma\left(\frac{d+1}{2\a}\right)}{2 \a(2\pi)^d} \abs x , & \forall x \in \R  
\end{disarray}
$$
where: $D_{\a}= \dis{\int_0^{\infty}} u^{2\a + \frac{d-1}{2}}2^{-(2\a + \frac{d-1}{2})} W_{0,\frac{d}{2}-1}(u) du$ and $W_{0,\nu}(z)= \frac{e^{-z/2}}{\Gamma(n+\frac{1}{2})} \int_0^{\infty} [t(1+\frac{t}{z})]^{\nu-\frac{1}{2}}e^{-t} dt$ is the Whittaker function.
\end{prop}
As in Proposition \ref{in1D}, the first inequality will be used to control $p_{\a}$ for large values, whereas the second one will be used to control $p_{\a}$ in the vicinity of $0$. The proof is based on \cite{Kolo} and \cite{sch}.

\begin{dem}
First we use the spherical coordinate system in dimension $d > 1$, that is to say we write: $\xi = (r, \theta, y)$ where $r$ belongs to the interval $(0,+\infty)$, $\theta$ belongs to $(0, \pi)$ and $\phi$ belongs to the $(d-2)$-sphere of radius $1$, with the main axis directed along $x$.
The jacobian of this change is 
 $J= r^{d-1} (\sin\theta)^{d-2} ,$ and we denote by $S_{d-2}$ the area of the $(d-2)$-sphere of radius $1$ ($\mathbb{S}_{d-2}$): $\mathbb{S}_{d-2}= 2 \dis{\frac{\pi^{\frac{d-1}{2}}}{\Gamma(\frac{d-1}{2})}}.$
Thus, we have, for $d>2$:
$$
\begin{disarray}{rcl}
p_{\a}(x) &=& \frac{1}{(2 \pi)^d} \int_0^{\infty} \int_{\theta=0}^{\pi} \int_{\phi \in \mathbb{S}_{d-2}} e^{-r^{2 \a}} e^{i\abs x r \cos\theta} r^{d-1} (\sin\theta)^{d-2} dr d\theta d\phi\\
&=& \frac{S_{d-2}}{(2 \pi)^d} \int_0^{\infty} \int_{-1}^1 e^{-r^{2 \a}} e^{i\abs x r t} r^{d-1} (1-t^2)^{\frac{d-3}{2}}dr dt\\
&=& \frac{S_{d-2}}{(2 \pi)^d} \int_0^{\infty} \int_{-1}^1 e^{-r^{2 \a}} \cos(\abs x r t) r^{d-1} (1-t^2)^{\frac{d-3}{2}}dr dt\\
\end{disarray}
$$
For $d=2$, we have:
$$
\begin{disarray}{rcl}
p_{\a}(x) &=&\frac{1}{(2 \pi)^2} \int_0^{\infty}\int_0^{2 \pi}  e^{-r^{2 \a}} e^{i\abs x r \cos\theta} rdr d\theta\\
&=&  \frac{2}{(2 \pi)^2} \int_0^{\infty}\int_0^{\pi}  e^{-r^{2 \a}} \cos(\abs x r \cos\theta) rdr d\theta\\
&=& \frac{S_{0}}{(2 \pi)^2} \int_0^{\infty} \int_{-1}^1 e^{-r^{2 \a}} \cos(\abs x r t) r (1-t^2)^{-\frac{1}{2}}dr dt\\
\end{disarray}
$$
and so we have the same expression.
Next, we use the Bessel function and the Whittaker function defined for $z \in \C \backslash \R_-$ and any real number $\nu > -\frac{1}{2}$ by the integral formulae:
$$ J_{\nu}(z)=\frac{(\nicefrac{z}{2})^{\nu}}{\Gamma(\nu+\frac{1}{2}) \sqrt \pi} \int_{-1}^1 (1-t^2)^{\nu-\frac{1}{2}} \cos(zt)dt, \quad W_{0,\nu}(z)= \frac{e^{-z/2}}{\Gamma(\nu+\frac{1}{2})} \int_0^{\infty} [t(1+\frac{t}{z})]^{\nu-\frac{1}{2}}e^{-t} dt.$$
Furthermore, for such a $\nu$ and  $z \in \R_+$, these functions are related by the formula:
$$J_{\nu}(z)= 2 \re \left(\frac{1}{\sqrt{2 \pi z}} e^{\frac{1}{2}( \nu + \frac{1}{2}) \pi i} W_{0, \nu}(2iz) \right).$$
Thanks to these special functions, and since $\frac{d}{2}-1 >-\frac{1}{2}$ , we can write:

\begin{eqnarray}\label{pa}
p_{\a}(x) &=&\frac{S_{d-2}}{(2 \pi)^d} \int_0^{\infty} e^{-r^{2 \a}} J_{\frac{d}{2}-1}(\abs x r) \Gamma \left(\frac{d-1}{2} \right) \sqrt{\pi} \frac{2^{\frac{d}{2}-1}}{\abs x^{\frac{d}{2}-1} r^{\frac{d}{2}-1} } r^{d-1}dr\nonumber\\
&=&(2 \pi)^{-\frac{d}{2}} \abs x^{1-{\frac{d}{2}} } \int_0^{\infty} e^{-r^{2 \a}} J_{\frac{d}{2}-1}(\abs x r) r^{\frac{d}{2}} dr \nonumber\\
&=&(2 \pi)^{-\frac{d}{2}} \abs x^{-d} \int_0^{\infty} e^{-\frac{y^{2 \a}}{\abs x^{2 \a}}} J_{\frac{d}{2}-1}(y) y^{\frac{d}{2}} dy \nonumber\\
&=&2 (2 \pi)^{-\frac{d}{2}} \abs x^{-d} \int_0^{\infty} e^{-\frac{y^{2 \a}}{\abs x^{2 \a}}} \re \left(\frac{1}{\sqrt{2 \pi y}} e^{\frac{1}{2}(\frac{d}{2}- \frac{1}{2}) \pi i} W_{0,\frac{d}{2}-1}(2iy) \right) y^{\frac{d}{2}} dy 
\end{eqnarray}
Now, we have to study: $I_{\a}(x) := \dis{\int_{0}^{\infty} }  e^{-\frac{y^{2 \a}}{\abs x^{2 \a}}}  e^{\frac{d-1}{4} \pi i} W_{0,\frac{d}{2}-1}(2iy) y^{\frac{d-1}{2}} dy $.
We follow the 1D method.
First, we rotate the line of integration by $-\frac{\pi}{4\a}$, to get:
$$
I_{\a}(x)= \int_0^{\infty}e^{i r^{2\a}\abs x^{-2\a}} r^{\frac{d-1}{2}}e^{-i \pi \frac{d+1}{8 \a}}e^{i \pi \frac{d-1}{4}}W_{0,\frac{d}{2}-1}(2ir e^{-i \frac{\pi}{4\a}} )dr.
$$
We denote by:
$$ I_{\a, \infty}(x)= \int_0^{\infty}(1+i r^{2\a}\abs x^{-2\a}) r^{\frac{d-1}{2}}e^{-i \pi \frac{d+1}{8 \a}}e^{i \pi \frac{d-1}{4}}W_{0,\frac{d}{2}-1}(2ir e^{-i \frac{\pi}{4\a}} )dr.$$
Recall that we have \eqref{pa}, so we have to take the real part of two terms:
\begin{itemize}
\item By rotating the integration line by  $-(2\a-1) \frac{\pi}{4\a}$ we get the real part of the first term :
$$
\begin{disarray}{rcl}
 I_{\a, \infty}^1(x)&:=& \int_0^{\infty} r^{\frac{d-1}{2}}e^{-i \pi \frac{d+1}{8 \a}}e^{i \pi \frac{d-1}{4}}W_{0,\frac{d}{2}-1}(2ir e^{-i \frac{\pi}{4\a}} )dr\\
&=&-i \int_0^{\infty} 2^{-\frac{d-1}{2}} u^{\frac{d-1}{2}} W_{0,\frac{d}{2}-1}(u) du\\
\end{disarray}
$$
If $u$ is a real number, $ W_{0,\frac{d}{2}-1}(u)$ is also a real number and consequently:
$\re \left( I_{\a,\infty}^1(x)\right)=0$.
\item The real part of the second term is computed with the same rotation:
$$
\begin{disarray}{rcl}
I_{\a, \infty}^2(x)&:=& \int_0^{\infty}i r^{2\a}\abs x^{-2\a}r^{\frac{d-1}{2}}e^{-i \pi \frac{d+1}{8 \a}}e^{i \pi \frac{d-1}{4}}W_{0,\frac{d}{2}-1}(2ir e^{-i \frac{\pi}{4\a}} )dr\\
&=&\int_0^{\infty} i e^{-i\a\pi} \frac{u^{2\a + \frac{d-1}{2}}}{2^{2\a + \frac{d-1}{2}}\abs x^{2\a}} W_{0,\frac{d}{2}-1}(u) du
\end{disarray}
$$
Consequently we have: $\re \left( I_{\a, \infty}^2(x) \right) = \sin(\a \pi) \dis{\int_0^{\infty}} \frac{u^{2\a + \frac{d-1}{2}}}{2^{2\a + \frac{d-1}{2}}\abs x^{2\a}} W_{0,\frac{d}{2}-1}(u) du.$
We denote by $D_{\a}$ the integral $\dis{\int_0^{\infty}} u^{2\a + \frac{d-1}{2}}2^{-(2\a + \frac{d-1}{2})} W_{0,\frac{d}{2}-1}(u) du$. From \cite{Erd}, we have:
 $D_{\a}=2^{2\a +\nicefrac{d}{2}-2} \Gamma(\a + \frac{d-1}{2} ) \Gamma(\a + \frac{1}{2}).$
\end{itemize}
Secondly, we write:
$$
I_{\a,r}(x)= \int_0^{\infty}\left(e^{i r^{2\a}\abs x^{-2\a}}-(1+i r^{2\a}\abs x^{-2\a}) \right) r^{\frac{d-1}{2}}e^{-i \pi \frac{d+1}{8 \a}}e^{i \pi \frac{d-1}{4}}W_{0,\frac{d}{2}-1}(2ir e^{-i \frac{\pi}{4\a}} )dr.
$$
Using the result obtained for $I_{\a,\infty}^1(x)$, we only have to treat the integral: (still denoted by $I_{\a,r}(x)$)
$$
I_{\a,r}(x)= \int_0^{\infty}\left(e^{i r^{2\a}\abs x^{-2\a}}-i r^{2\a}\abs x^{-2\a} \right) r^{\frac{d-1}{2}}e^{-i \pi \frac{d+1}{8 \a}}e^{i \pi \frac{d-1}{4}}W_{0,\frac{d}{2}-1}(2ir e^{-i \frac{\pi}{4\a}} )dr.
$$
Introducing the new variable $u=r^{2\a}$:
$$
I_{\a,r}(x)=(2\a)^{-1} \int_0^{\infty}\left(e^{i u\abs x^{-2\a}}-i u\abs x^{-2\a} \right) u^{\frac{d+1}{4\a}-1}e^{-i \pi \frac{d+1}{8 \a}}e^{i \pi \frac{d-1}{4}}W_{0,\frac{d}{2}-1}(2iu^{\nicefrac{1}{2\a}} e^{-i \frac{\pi}{4\a}} )du.
$$
Then, we introduce: $f_u(\a)=u^{\frac{d+1}{4\a}-1} e^{-i \pi \frac{d+1}{8 \a}}W_{0,\frac{d}{2}-1}(2iu^{\nicefrac{1}{2\a}} e^{-i \frac{\pi}{4\a}} ),$ and write:
$$
\begin{disarray}{rcl}
I_{\a,r}(x)&=&(2\a)^{-1} \int_0^{\infty}e^{i u\abs x^{-2\a}} e^{i \pi \frac{d-1}{4}} f_u(1)du - (2\a)^{-1} \int_0^{\infty}i u\abs x^{-2\a} e^{i \pi \frac{d-1}{4}} f_u(1)du \\
&&+ (2\a)^{-1} \int_0^{\infty}\left(e^{i u\abs x^{-2\a}}-i u\abs x^{-2\a} \right)(f_u(\a)-f_u(1))e^{i \pi \frac{d-1}{4}}du.
\end{disarray}$$
Let us take the real part of each element of the right hand side:
\begin{itemize}
\item  $(2\a)^{-1}\dis{\int_0^{\infty}}e^{i u\abs x^{-2\a}} e^{i \pi \frac{d-1}{4}} f_u(1)du$ is treated by introducing $y= u^{\frac{1}{2}}e^{-\frac{i\pi}{4}}$:
$$
\begin{disarray}{rcl}
\re \left((2\a)^{-1}\dis{\int_0^{\infty}}e^{i u\abs x^{-2\a}} e^{i \pi \frac{d-1}{4}} f_u(1)du \right)
&= & \a^{-1}\re \left( \int_0^{\infty} e^{-\frac{y^{2 }}{\abs x^{2 \a}}}  e^{\frac{d-1}{4} \pi i} W_{0,\frac{d}{2}-1}(2iy) y^{\frac{d-1}{2}} dy \right)\\
&=& \a^{-1} p_{1}(x^{\a})2^{-1} (2 \pi)^{\frac{d+1}{2}} \abs x^{d}\\
&=& \a^{-1} e^{-\frac{\abs x ^{2\a}}{4}} \abs x^{d \a} \sqrt \pi 2^{-\frac{d+1}{2}}\\
\end{disarray}
$$
\item $(2\a)^{-1} \dis{ \int_0^{\infty}}i u\abs x^{-2\a} e^{i \pi \frac{d-1}{4}} f_u(1)du$ is treated by introducing: $u=-ir^2$. Thus we obtain:
$$
(2\a)^{-1} \dis{ \int_0^{\infty}}i u\abs x^{-2\a} e^{i \pi \frac{d-1}{4}} f_u(1)du=\a^{-1} \int_0^{\infty}i \abs x^{-2\a} r^{\frac{d+1}{2}} W_{0,\frac{d}{2}-1}(2r) r dr
$$
Consequently: $\re \left( (2\a)^{-1} \dis{ \int_0^{\infty}}i u\abs x^{-2\a} e^{i \pi \frac{d-1}{4}} f_u(1)du \right) = 0$.
\item We write:
$$
\begin{array}{l}
(2\a)^{-1} \dis{\int_0^{\infty}}\left(e^{i u\abs x^{-2\a}}-i u\abs x^{-2\a} \right)(f_u(\a)-f_u(1))e^{i \pi \frac{d-1}{4}}du =  \\
 \hspace{4cm} (2\a)^{-1} \dis{\int_0^{\infty}}\left(e^{i u\abs x^{-2\a}}-1-i u\abs x^{-2\a} \right)(f_u(\a)-f_u(1))e^{i \pi \frac{d-1}{4}}du.
\end{array}
$$

 Then this integral is negligible since it is bounded by  $C(1-\a) \abs x^{-4\a}$, where $C$ is a constant independent of $\a$. Indeed:
$$
\begin{disarray}{l}
\abs{\re \left( (2\a)^{-1} \dis{\int_0^{\infty}}\left(e^{i u\abs x^{-2\a}}-1-i u\abs x^{-2\a} \right)(f_u(\a)-f_u(1))e^{i \pi \frac{d-1}{4}}du\right)}\\
 \hspace{7cm}\leq C\dis{\int_0^{\infty}}u^2 \abs x^{-4\a}(1-\a) \sup_{y\in (\a,1)} \abs{\partial_y f_u(y)}du\\
\end{disarray}
$$
Moreover, for $y\in (\a, 1)$:
$$
\begin{disarray}{rcl}
\partial_y f_u(y)&=&-\frac{1}{8y^2} \left( iW_{1,\frac{d}{2}-1}(2iu^{\nicefrac{1}{2y}}e^{-i\frac{\pi}{4y}})-i(d+1)W_{0,\frac{d}{2}-1}(2iu^{\nicefrac{1}{2y}}e^{-i\frac{\pi}{4y}})\right.\\
&&\left.+2 W_{0,\frac{d}{2}-1}(2iu^{\nicefrac{1}{2y}}e^{-i\frac{\pi}{4y}})u^{\nicefrac{1}{2y}}e^{-i\frac{\pi}{4y}} \right) u^{\frac{d+1}{4y}-1} \left(-i\pi +2\ln(u) \right)e^{-i\pi \frac{d+1}{8y}}.\\
\end{disarray}
$$
So:
$$
\begin{disarray}{rcl}
\abs{\partial_y f_u(y)}& \leq &\frac{1}{8\a^2}\left( \abs{W_{1,\frac{d}{2}-1}(2iu^{\nicefrac{1}{2y}}e^{-i\frac{\pi}{4y}})}+\abs{W_{0,\frac{d}{2}-1}(2iu^{\nicefrac{1}{2y}}e^{-i\frac{\pi}{4y}})} (1+d+2u^{\frac{d+1}{4y}-1}) \right)\\
&&(\pi +2 \abs{\ln(u)}).
\end{disarray}
$$
And thus (see the 1D study):
$$ \dis{\int_0^{\infty}}u^2 \abs x^{-4\a}(1-\a) \sup_{y\in (\a,1)} \abs{\partial_y f_u(y)}du
\leq C(1-\a) \abs x^{-4\a}.$$
\end{itemize}
Consequently:
$$
\begin{disarray}{rcl}
\re(I_{\a}(x))&=&\sin(\a \pi) D_{\a} \abs x^{-2\a}+\a^{-1} e^{-\frac{\abs x ^{2\a}}{4}} \abs x^{d \a} \sqrt \pi 2^{-\frac{d+1}{2}}\\
&& +\re\left((2\a)^{-1} \dis{\int_0^{\infty}}\left(e^{i u\abs x^{-2\a}}-i u\abs x^{-2\a} \right)(f_u(\a)-f_u(1))e^{i \pi \frac{d-1}{4}}du\right).
\end{disarray}
$$
Recall that we have:
$$p_{\a}(x)=2 (2 \pi)^{-\frac{d+1}{2}} \abs x^{-d}  \re \left( I_{_a}(x) \right).$$
So, we obtain the inequality:
$$\abs{p_{\a}(x) - 2(2\pi)^{-\frac{d+1}{2}}\sin(\a \pi) D_{\a}\abs x^{-(d+2\a)}- (4 \pi)^{-\frac{d}{2}} \abs x^{-(1-\a)d}\a^{-1} e^{-\frac{\abs x ^{2\a}}{4}}}\leq C(1-\a)\abs x^{-(d+4\a)} $$
To get an estimate for $x$ in a compact set, we use the expression of $p_{\a}(0)$ given by:
$$
\begin{disarray}{rcl}
p_{\a}(0) =  \frac{1}{(2 \pi)^d}\int_{\R^d} e^{-\abs \xi^{2 \a}} d\xi = \frac{\Gamma(\nicefrac{d}{2\a}) }{2 \a (2\pi)^d} < {+\infty}
\end{disarray}
$$
For all $y \in \R$, we have:
$$
\begin{disarray}{rclcl}
 \abs{p_{\a}'(y)} & \leq& \frac{1}{(2\pi)^d} \int_{\R^d} \abs {\xi} e^{-\abs{\xi}^{2\a}}d\xi &=&\frac{1}{2 \a(2\pi)^d} \Gamma\left(\frac{d+1}{2\a}\right)\\
\end{disarray}
$$
So: 
$$
\begin{disarray}{rclr}
  \abs{p_{\a}(x)-p_{\a}(0)} &\leq&  \frac{\Gamma\left(\frac{d+1}{2\a}\right)}{2 \a(2\pi)^d} \abs x, & \forall x \in \R. 
\end{disarray}
$$

\end{dem}
\subsection{Consequence: Heuristics for the level set $\{ x\in \R^d, \ u(x,t)=\frac{1}{2}\}$ for $t$ large enough}
The study of the level set relies on the fact that the function: $y \mapsto y^{(d+2)\a} e^{-\frac{y^{2\a}}{4}}$ is non increasing for $y \geq [2(d+2)]^{\nicefrac{1}{(2\a)}}$ (value for which the maximum is reached). We denote by $\xi_{\a}$ the solution to: $y^{(d+2)\a} e^{-\frac{y^{2\a}}{4}}= C_{\a}\sin(\a \pi) $ larger than $[2(d+2)]^{\nicefrac{1}{(2\a)}}$, where $C_{\a}=\frac{2}{\sqrt{\pi}}$ if $d=1$, and $C_{\a}=\frac{2^{\frac{d+1}{2}} D_{ \a}\a}{\sqrt{\pi}}$ if $d>1$.
Define: $$\tau_{\a}= \dis{\frac{\xi_{\a}^{2\a}}{4}}\underset{\a \rightarrow 1}{\sim}- \ln(1-\a).$$
Set $\tild{C}^{2\a}\geq\dis{\frac{C (1-\a)}{\sin(\a \pi)}}$ a constant independent of $\a$. We notice that: 
\begin{itemize}
\item For $\tild{C} \leq \abs{\xi} \leq \xi_{\a}$, we have:
$$\frac{e^{-\frac{\abs{\xi}}{4}^{2\a}}}{\abs{\xi}^{d(1-\a)}} \geq C\frac{\sin(\a \pi) }{ \abs{\xi}^{d+2\a}} \quad \mbox { and } \quad C\frac{1-\a}{\abs{\xi}^{d+4\a}} \leq  \frac{\sin(\a \pi) }{ \abs{\xi}^{d+2\a}}.$$
So the fundamental solution will behave like: $ \dis{\frac{e^{-\frac{\abs{x}}{4t}^{2\a}}}{t^{\nicefrac{d}{2}}\abs x^{d(1-\a)}}} $, which is very close to the heat kernel when $\a$ tends to $1$: the propagation should be linear (see \cite{AW}).
\item There exists $\a_1 \in (\frac{1}{2},1)$ such that $\forall \a \in (\a_1,1)$, for $\abs{\xi} \geq \xi_{\a}$,  we have:
$$\frac{e^{-\frac{\abs{\xi}}{4}^{2\a}}}{\abs{\xi}^{d(1-\a)}} \leq C \frac{\sin(\a \pi)}{ \abs{\xi}^{d+2\a}} \quad \mbox { and } \quad
 C\frac{1-\a}{\abs{\xi}^{d+4\a}} \leq  \frac{\sin(\a \pi) }{ \abs{\xi}^{d+2\a}}.$$
So the fundamental solution will behave like: $\dis{\frac{\sin(\a \pi) t}{ \abs{x}^{d+2\a}}}, $ and as is shown in \cite{JMRXC} the propogation should be exponential.
\end{itemize}
\section{An intermediate result}\label{inter}
In the forthcoming parts \ref{propexp1} and \ref{propexp}, we will study the evolution of the level set $\{ x \in \R^d \  | \ u(x,t)=\eps_{\a} \}$ with $\eps_{\a}=\sin(\a \pi)^{1+ \kappa}, \kappa >0$. However, we need to know how the level set $\{ x \in \R^d \  | \ u(x,t)=\underline{\eps} \}$, $\underline{\eps} >0$ small and independent of $\a$, evolves. The following lemma makes the connection between these two level sets.
Let us consider the evolution problem:
\begin{equation} \label{v}
\left\{
\begin{array}{lclr} 
v_t+(-\Delta )^{\a}v&=&v-v^2,& x \in B_M , t >0 \\
v(x,t)&=&0, & x \in \R^d \setminus B_M, t >0\\
v(x, 0)&=&\eps_{\a} \mathds{1}_{B_{M-1}}(x), & x \in \R^d\\
\end{array}
\right.
\end{equation}
for $\eps_{\a}= \sin(\a \pi)^{1+\kappa}$, for $\kappa>0$ and $M$ large enough so that the principal Dirichlet eigenvalue of $(-\Delta )^{\a}-I$ in $B_M$ is negative. This is possible due to Theorem 1.1 in \cite{BRR}: the first eigenvalue of $(-\Delta)^{\a}-I$ with Dirichlet condition outside $B_M$ tends to $-1$ as $M$ tends to $+\infty$.
\begin{prop}\label{youpi}
There exist a constant $c>0$ independent of $\a$,  $\tild{\tau_{\a}}<c \tau_{\a}$, $\underline{\eps} \in (0,1)$ independent of $\a$, and $m \in B_M$ such that:
$$v(m, \tild{\tau_{\a}}) \geq \underline{\eps}.$$
\end{prop}
\begin{dem}
We have: $v_t +A_{\a}v=-v^2$, where $A_{\a}=(-\Delta )^{\a}-I$ is self-adjoint and its principal eigenvalue is denoted by $\mu_1^{\a}<0$.
 In $L^2(B_M)$, let $e_1^{\a}$ be an element of an eigenvector basis,  corresponding to $\mu_1^{\a}$.
Since $e_1^{\a}>0$ in $B_M$, there exists $C_0>0$ such that:
$$v(x,0) \leq C_0 \eps_{\a} e_1^{\a}(x), \mbox{ for } x \in B_{M-1}.$$
Thus, $\overline{v}(x,t)= C_0 \eps_{\a} e_1^{\a}(x)e^{-\mu_1^{\a}t}$ is a super solution to \eqref{v} in $B_M$, coincides with $v$ outside, and $\norme{\overline{v}(\cdot,t)}_2=C_0 \eps_{\a}e^{-\mu_1^{\a}t}, \forall t >0$.
 Let $\underline{\eps} \in (0,1)$ independent of $\a$. The norm $\norme{\overline{v}(\cdot,t)}_2$ reaches $2 \sqrt{|B_{M}|}\underline{\eps}$ for $t=\tild{\tau_{\a}}= (-\mu_1^{\a})^{-1} \ln(\frac{2 \sqrt{|B_{M}|}\underline{\eps}}{B\eps_{\a}})$ (this computation will be done later, see $\eqref{T_0}$). Using the expression of $\eps_{\a}$, we have the existence of a constant $c$ independent of $\a$ such that:  $\tild{\tau_{\a}}<c \tau_{\a}$. Then, define $w$ by: $w(x,t)=\overline{v}(x,t)-v(x,t)$, solution to:
$$
\left\{
\begin{array}{lclr}
w_t+(-\Delta )^{\a}w&=&v^2,& x \in B_{M}, t >0\\
w(x,t)&=&0, &  x \in \R^d \setminus B_M, t >0\\
w(x, 0)&=&C_0 \eps_{\a} e_1^{\a}(x)- \eps_{\a} \mathds{1}_{B_{M-1} }(x), & x \in \R^d\\
\end{array}
\right.
$$
Let $\widehat{\tau_{\a}}$ be the largest time for which the following holds:
$$\norme{w(\cdot,t)}_2 \leq \frac{\norme{\overline{v}(\cdot,t)}_2}{2}, \ \forall t \leq \widehat{\tau_{\a}}. $$
Assume $\widehat{\tau_{\a}}< \tild{\tau_{\a}}$. In the following inequalities, we use the fact $e_1^{\a}$ is continuous and its $L^2$ and $L^{\infty}$ norms are comparable. Moreover, since $A_{\a}$ is symmetric, we have :
$$\norme{ e^{-A_{\a}t}}_{L^2 \rightarrow L^2} \leq e^{-\mu_1^{\a}t}, \ \forall t>0.$$
Thus, for $t \leq \widehat{\tau_{\a}}$ and $\a$ close to $1$ enough:

$$
\begin{disarray}{rcl}
\norme{w(\cdot,t)}_2&=&  \norme{\int_0^t e^{-A_{\a}(t-s)}v(\cdot,s)^2ds+e^{-A_{\a}t}\left(C_0 \eps_{\a} e_1^{\a}(x)- \eps_{\a} \mathds{1}_{B_{M-1} }(x)\right)}_2\\
&\leq&\int_0^t \norme{ e^{-A_{\a}(t-s)}}_{L^2 \rightarrow L^2} \norme{v(\cdot,s) \overline{v}(\cdot,s)}_2ds +C \eps_{\a} \norme{ e^{-A_{\a}t}}_{L^2 \rightarrow L^2}  \\
&\leq&C\int_0^t \norme{ e^{-A_{\a}(t-s)}}_{L^2 \rightarrow L^2} (\norme{\overline{v}(\cdot,s)}_{\infty}^2+\norme{ \overline{v}(\cdot,s)}_{\infty}\norme{w(\cdot,s)}_2)ds+C \eps_{\a} e^{-\mu_1^{\a}t}\\
&\leq&C\int_0^t e^{-\mu_1^{\a}(t-s)} (\norme{\overline{v}(\cdot,s)}_{2}^2+\norme{ \overline{v}(\cdot,s)}_{2}\norme{w(\cdot,s)}_2)ds+C \eps_{\a} e^{-\mu_1^{\a}t}\\
&\leq&C\int_0^t e^{-\mu_1^{\a}(t-s)}\norme{ \overline{v}(\cdot,s)}_{2}^2ds+C \eps_{\a} e^{-\mu_1^{\a}t}\\
&\leq&C\int_0^t e^{-\mu_1^{\a}(t-s)}C_0^2 \eps_{\a}^2e^{-2\mu_1^{\a}s}ds+C \eps_{\a} e^{-\mu_1^{\a}t}\\
&\leq&CC_0^2\eps_{\a}^2\frac{e^{-2\mu_1^{\a}t}}{\abs{\mu_1^{\a}}}+C \eps_{\a} e^{-\mu_1^{\a}t}\\
&\leq&C\eps_{\a}e^{-\mu_1^{\a}t}(\norme{ \overline{v}(\cdot,t)}_{2}+1),
\end{disarray}
$$
For $t=\widehat{\tau_{\a}}$, using the fact $\mu_1^{\a}<0$ and $\widehat{\tau_{\a}}< \tild{\tau_{\a}}$:
$$\norme{w(\cdot,\widehat{\tau_{\a}})}_2=\frac{\norme{\overline{v}(\cdot,\widehat{\tau_{\a}})}_2}{2} \leq 2\sqrt{|B_{M}|}C\frac{\underline{\eps}}{B}(\norme{ \overline{v}(\cdot,\widehat{\tau_{\a}})}_{2}+1),$$
taking $\underline{\eps}$ smaller if necessary, we get a contradiction. Consequently  $\widehat{\tau_{\a}}\geq \tild{\tau_{\a}}$, and so:
$$\norme{w(\cdot,\tild{\tau_{\a}})}_2 \leq \frac{\norme{\overline{v}(\cdot,\tild{\tau_{\a}})}_2}{2} .$$
Thus:
$$\norme{v(\cdot,\tild{\tau_{\a}})}_2 \geq \norme{\overline{v}(\cdot,\tild{\tau_{\a}})}_2-\norme{w(\cdot,\tild{\tau_{\a}})}_2 \geq \frac{\norme{\overline{v}(\cdot,\tild{\tau_{\a}})}_2}{2}.$$
However: $\norme{v(\cdot,\tild{\tau_{\a}})}_2\leq \sqrt{|B_{M}|}\norme{v(\cdot,\tild{\tau_{\a}})}_{\infty},$ so there exists $m\in B_M$ such that: $$v(m,\tild{\tau_{\a}})\geq \underline{\eps}.$$
\end{dem}

\section{Initial data with compact support} \label{4}
This section contains the proof of Theorem \ref{thm}. First we are concerned with the linear propagation phase and then the exponential one. The difficulty is to find a lower bound of the solution. The idea is to use an iterative scheme and in each iteration we truncate the solution so that the reaction term in \eqref{sys} is bigger than a linear term. Then, the study of the fundamental solution in Section \ref{2} leads to the result.

\subsection{The linear propagation phase}\label{lin}
We consider the evolution problem:\\
\begin{equation}\label{systeme}
\left\{
\begin{array}{rcl}
u_t +(-\Delta)^{\a}u&=&u-u^2, \quad \R^d, t>0\\
u(x,0)&=&u_0(x), \hfill x \in \R^d\\
\end{array}
\right.
\end{equation}
where $u_0$ is compactly supported, continuous and $u_0 \in [0,1].$ The following lemma and its corollary are inspired by Cabr\'e and Roquejoffre in \cite{JMRXC}.
\begin{lemme}\label{it}
For every $0<\sigma<2$ and $\a \in (\nicefrac{1}{2},1)$, there exist $\eps_0 \in (0,1)$ and $T_0 \geq 1$ depending only on $\sigma$ and $\eps_0$ for which the following holds.
Given $r_0 \in (1,C\tau_{\a})$, $C$ independent of $\a$ and $\eps \in (0, \eps_0)$, let $\underline{u_0}=\eps \mathds{1}_{B_{r_0}(0)}$.
Then, the solution to \eqref{systeme} with initial condition $\underline{u_0}$ satisfies, for all $k\in \N$ such that $kT_0<\tau_{\a}$:
$$u(x,kT_0) \geq \eps \mbox{  for  } \abs x \leq r_0 + k\sigma T_0^{\nicefrac{1}{\a}} .$$
\end{lemme}
\begin{dem}
For $k=0$, the result is obvious.\\
For $k=1$, we notice: $u-u^2 \leq u$, for $u\in [0,1]$; moreover for every $\delta \in (0,1)$, as long as $u\leq \delta$ we have : $u-u^2 \geq (1-\delta)u$.Thus, taking $\delta >> \eps$ we have a super solution and a sub solution to \eqref{systeme}:
$$ \underline{u}(x,t):=e^{(1-\delta)t} \int_{\R^d} \underline{u_0}(y) p(x-y,t) dy \leq u(x,t) \leq e^{t} \int_{\R^d} \underline{u_0}(y) p(x-y,t) dy=:\overline{u}(x,t),$$
as long as $\overline{u} \leq \delta.$
Let $\delta >> \eps$, $T_0 > 0$ chosen so that:
\begin{equation}\label{T}
\forall t\in (0,T_0), \forall x \in \R^d, \quad \overline{u}(x,t) \leq \delta.
\end{equation}
Thanks to Lemma 2.3 in \cite{JMRXC},  we have : If $u$ and $v$ : $\R^d \rightarrow \R$, $u \in L^1, v \in L^{\infty}$ are positive radially symmetric and nonincreasing functions, then $u \star v $ is also  positive radially symmetric and nonincreasing. Here, thanks to Proposition \ref{inD} applied to $p(x,t)=t^{-\frac{d}{2\a}}p_{\a}(x t^{-\frac{1}{2\a}})$, we have that $p$ is smaller than the function: 
$$x \longmapsto 
\left\{\begin{array}{lr}
\dis{C \left(\frac{(1-\a)t^2}{\abs x ^{d+4\a}}+\frac{\sin(\a\pi)t}{\abs x ^{d+2\a}}+ \frac{e^{- \frac{\abs x}{4t}^{2\a}}}{ t^{\nicefrac{d}{2}}\abs x^{d(1-\a)}}\right)},& \mbox{ if } \abs x \geq \tild{C} t^{\nicefrac{1}{(2\a)}}\\
\hat{C} t^{-\nicefrac{d}{2\a}}, & \mbox{ if } \abs x \leq \tild{C} t^{\nicefrac{1}{(2\a)}},
\end{array} \right.
$$
where $\hat{C}$ is a constant large enough and independent of $\a$. This function is a positive, radially symmetric and nonincreasing function.
So, with the initial condition $\underline{u_0}=\eps  \mathds{1}_{B_{r_0}(0)}$, it is sufficient to estimate $\overline{u}(0,t)$.
As a consequence, we have to find $T_0$ such that : $\forall t\in (0,T_0), \quad \overline{u}(0,t) \leq \delta$, to get \eqref{T}. Yet, the inequalities on $p$ lead to:
\begin{eqnarray}\label{T_0}
\overline{u}(0,t) &\leq& e^t \int_{\R^d} \eps \mathds{1}_{B_{r_0}(0)}(y) p(-y,t)dy \nonumber \\
&\leq& \eps e^t \int_{\abs y \leq \tild{C} t^{\nicefrac{1}{2\a}}}\hat{C} t^{-\nicefrac{d}{2\a}}dy+ C \eps e^t \int_{\abs y \geq \tild{C} t^{\nicefrac{1}{2\a}}} \frac{(1-\a)t^2}{\abs y ^{d+4\a}}+\frac{ \sin(\a\pi)t}{\abs y ^{d+2\a}}\nonumber + \frac{e^{-\frac{\abs y}{4t}^{2\a}}}{t^{\nicefrac {d}{2}}\abs y^{d(1-\a)}}dy \nonumber\\
&\leq&  C  \eps e^t \left(1+  \frac{(1-\a)}{2 \a \pi \tild{C} ^{4\a}} +  \frac{\sin(\a\pi)}{\a\tild{C} ^{2\a}} +\int_{\abs z \geq \tild{C}2^{-\nicefrac{1}{\a}}} \frac{ e^{-\abs z^{2\a}}}{\abs z^{d(1-\a)}} dz \right)\nonumber\\
&\leq & C \eps e^t ,
\end{eqnarray}
where $C$ is independent of $\a$. Thus, $T_0 = \ln\left( \dis{\frac{\delta}{B\eps}} \right)$ is smaller than $ \tau_{\a}= - \ln(1-\a)$  (taking $1-\a$ smaller if necessary) and we have: 
$$\forall t\in (0,T_0), \forall x \in \R^d, \quad u(x,t) \leq \delta.$$
We notice that the smaller  $\eps$ is, the larger  $T_0$  is, but we always have \eqref{T}. 

Then, we look for $r_1> r_0$ for which we have: $u(x,T_0) \geq \eps,  \forall \abs x \leq r_1$. 
To find it, we look for $x_1>r_0$ so  that $\underline{u}$ is larger than $\eps$ for $\abs x$ smaller than $x_1$.
First, we notice that, for $y \in \R^d$, if $\abs{y-x} \geq \tild{C} T_0^{\nicefrac{1}{2\a}}$, where $\tild{C}^{2\a}\geq\dis{\frac{C (1-\a)}{\sin(\a \pi)}}$, $\tild{C}$ independent of $\a$, then: 
\begin{equation}\label{neg}
 \frac{C(1-\a) T_0^2}{\abs{y-x}^{d+4\a}} \leq \frac{\sin(\a \pi) T_0}{\abs{y-x}^{d+2\a}}.
\end{equation}
Next, we prove the existence of a constant $\overline{C} > 4^{\nicefrac{1}{2\a}}$ so that:
$$\tild{C}T_0^{\nicefrac{1}{2\a}} \leq  \abs{x_1-r_0} \leq \overline{C} T_0^{\nicefrac{1}{\a}}.$$
Indeed:
\begin{itemize}
\item if for all $\overline{C} >  4^{\nicefrac{1}{2\a}}$, $\abs{x_1-r_0} > \overline{C} T_0^{\nicefrac{1}{\a}}$, then $\overline{u}(x,T_0)$ is strictly smaller than $\eps$, which is impossible: indeed, denoting by $e_1$ the first vector of the standard basis of $\R^d$, for $y \in \R^d$ such that $ \abs y \leq r_0$, we get $\abs{x_1e_1-y}\geq \abs{x_1-r_0} > \overline{C} T_0^{\nicefrac{1}{\a}}  $ and for $\tild{x_1}=x_1e_1$:
$$
\begin{disarray}{rcl}
\overline{u}(\tild{x_1},T_0) &\leq& C e^{T_0}\eps \int_{\abs y \leq r_0} \frac{{T_0}^2 (1-\a) }{\abs{\tild{x_1}-y}^{d+4\a}}dy + C e^{T_0}\eps  \int_{\abs y \leq r_0} \frac{T_0 \sin (\a \pi)}{\abs{\tild{x_1}-y}^{d+2\a}} dy \\
&& + C e^{T_0}\eps\int_{\abs y \leq r_0}\frac{e^{-\frac{\abs{y-\tild{x_1}}}{4T_0}^{2\a}}}{T_0^{\nicefrac{d}{2}}\abs{\tild{x_1}-y}^{d(1-\a)}} dy \\
&\leq&   \frac{C \eps {T_0}^2}{4\a \abs{x_1-r_0}^{4\a}} +   \frac{C\eps T_0  }{2\a \abs{x_1-r_0}^{2\a}}  +C\eps  e^{T_0} \int_{\abs{z} \geq \frac{x_1-r_0}{(4T_0)^{\nicefrac{1}{2\a}}}} e^{- \abs{z}^{2\a}}\abs{z}^{d(\a-1)}dz \\
&\leq&  C \eps  T_0^{-2} +  C \eps T_0^{-1}  +C\eps  e^{T_0-\frac{\overline{C}^{2\a}}{4} T_0}  \\
&<& \eps, 
\end{disarray}
$$
for $T_0$ large enough, since $e^{T_0} \leq C\sin(\a \pi)^{-1} $. Here we use the fact that for all $y$ large enough, there exists a constant $C$ independent of $\a$ such that: 
\begin{equation}\label{equivalent}
\int_{\abs z>y}e^{- \abs z^{2\a}}\abs z^{d(\a-1)}dz \leq C \frac{e^{-y^{2\a}}}{2\a y^{\a(2-d)}}. 
\end{equation}
 So, there exists $\overline{C} >  4^{\nicefrac{1}{2\a}}$ such that $\abs{x_1-r_0} \leq \overline{C} T_0^{\nicefrac{1}{\a}}$.
\item if $ \abs{x_1-r_0} \leq \tild{C}T_0^{\nicefrac{1}{2\a}}$, then $\underline{u}$ is larger than $2 \eps$, which is impossible: indeed, we have  $\abs{(x_1- \tild{C}T_0^{\nicefrac{1}{2\a}})e_1} \leq r_0$, so: 
$$
\begin{disarray}{rcl}
\underline{u}(x_1e_1,T_0) &=& e^{(1-\delta)T_0} \int_{\R^d} \underline{u_0}(y) p(x_1e_1-y,T_0) dy\\
&\geq& e^{(1-\delta)T_0} \eps  p( \tild{C}T_0^{\nicefrac{1}{2\a}}e_1,T_0) \\
&\geq&C e^{(1-\delta)T_0} \eps e^{-\frac{\tild{C}^{2\a}}{4}}T_0^{-\nicefrac{d}{2}}\\
&\geq& 2 \eps, \\
\end{disarray}
$$
the last inequality is obtained taking $T_0$ larger if necessary with $T_0 < \tau_{_a}$.
\end{itemize}
Thus, using the fact the heat kernel is positive and Proposition \ref{inD}, we have for all $x$ in $\R^d$ and for $r_0 > \tild{C} T_0^{\nicefrac{1}{2\a}}$:
$$
\begin{disarray}{rcl}
u(x,T_0) &\geq& e^{(1-\delta)T_0} \int_{\R^d} p(x-y,T_0) \underline{u_0}(y) dy \\
&\geq&e^{(1-\delta)T_0}\eps\int_{\tiny{\begin{array}{l}\abs{y-x}\leq \tild{C} T_0^{\nicefrac{1}{2\a}} \\ \abs y \leq r_0\end{array}}}\left(\frac{p_{\a}(0)}{T_0^{\nicefrac{d}{2\a}}}-\frac{C\abs{x-y} }{T_0^{\nicefrac{(d+1)}{2\a}}}\right)dy\\ &&+Ce^{(1-\delta)T_0}\eps\int_{\tiny{\begin{array}{l}\abs{y-x}\geq \tild{C} T_0^{\nicefrac{1}{2\a}} \\ \abs y \leq r_0\end{array}}}\frac{e^{-\frac{\abs{y-x}}{4T_0}^{2\a}}}{ T_0^{\nicefrac{d}{2}}\abs{y-x}^{d(1-\a)}} dy\\
&\geq&e^{(1-\delta)T_0}\eps\int_{\tiny{\begin{array}{l}\abs{y-x}\leq \tild{C} T_0^{\nicefrac{1}{2\a}} \\ \abs y \leq r_0\end{array}}}\hat{C}T_0^{-\nicefrac{d}{2\a}}dy\\ &&+Ce^{(1-\delta)T_0}\eps\int_{\tiny{\begin{array}{l}\abs{y-x}\geq \tild{C} T_0^{\nicefrac{1}{2\a}} \\ \abs y \leq r_0\end{array}}}\frac{e^{-\frac{\abs{y-x}}{4T_0}^{2\a}}}{T_0^{\nicefrac{d}{2}}\abs{y-x}^{d(1-\a)}} dy=:w(x,T_0).\\
\end{disarray}
$$
Notice here that $\hat{C}>0$. Let us now study $w$ since it is radially symmetric and nonincreasing. We have for $x \in \R^d$ such that:  $\tild{C} T_0^{\nicefrac{1}{2\a}} \leq \abs{x-r_0e_1} \leq \overline{C} T_0^{\nicefrac{1}{\a}} $:
$$
\begin{disarray}{rcl}
w(x,T_0) &\geq&Ce^{(1-\delta)T_0}\eps\int_{\tiny{\begin{array}{l}\abs{y-x}\geq \tild{C} T_0^{\nicefrac{1}{2\a}} \\ \abs y \leq r_0\end{array}}}\frac{e^{-\frac{\abs{y-x}}{4T_0}^{2\a}}}{T_0^{\nicefrac{d}{2}}\abs{y-x}^{d(1-\a)}} dy\\
&\geq&Ce^{(1-\delta)T_0}\eps\frac{e^{-\frac{\abs{r_0e_1-x}}{4T_0}^{2\a}}}{T_0^{\nicefrac{d}{2}}\abs{r_0e_1-x}^{d(1-\a)}} \\
&\geq& C\eps\frac{e^{(1-\delta)T_0-\frac{\abs{r_0e_1-x}}{4T_0}^{2\a}}}{ T_0^{\nicefrac{d}{\a}-\nicefrac{d}{2}}}.
\end{disarray}
$$
Let us define $x_1$ by: $$C\eps\frac{e^{(1-\delta)T_0-\frac{\abs{r_0-x_1}^{2\a}}{4T_0}}}{ T_0^{\nicefrac{d}{\a}-\nicefrac{d}{2}}}=\eps,$$
that is to say: $$ x_1= r_0+2^{\nicefrac{1}{\a}}T_0 ^{\nicefrac{1}{\a}}\left(1-\delta-\frac{1}{4T_0}\ln \left(CT_0^{\nicefrac{d}{\a}-\nicefrac{d}{2}}\right) \right)^{\nicefrac{1}{2\a}}.$$
Consequently, for $\sigma < 2$, we take $\delta$ small enough, $T_0$ large enough (but smaller than $\tau_{\a}$) and $\a$ close enough to $1$ so that:
$$\sigma < 2^{\nicefrac{1}{\a}} \left(1-\delta-\frac{1}{4T_0}\ln \left(CT_0^{\nicefrac{d}{\a}-\nicefrac{d}{2}}\right) \right)^{\nicefrac{1}{2\a}}.$$
Now, let us define $r_1:=r_0+\sigma T_0^{\nicefrac{1}{\a}}>r_0$ so that $ r_1<x_1 $ and since $w$ is radially symmetric and nonincreasing:
$$u(r_1e_1,T_0) \geq w(r_1e_1,T_0) \geq w(x_1e_1,T_0) = \eps. $$
And:
$$ u(x,T_0) \geq w(x,T_0) \geq w(r_1e_1,T_0) \geq \eps, \ \forall \abs x \leq r_1.$$
Finally, $u(.,T_0) \geq \eps \mathds{1}_{B_{r_1}(0)}=: \underline{\underline{u_0}}$ and we can repeat the argument above, now with initial time $T_0$ and inital condition $\underline{\underline{u_0}}$ as long as $kT_0 < \tau_{\a}$ and get that:
$$u(x, kT_0)\geq \eps \mbox{ for } \abs x \leq r_k,$$
for all $k \in \N$ satisfying $kT_0 < \tau_{\a}$, with: $r_k \geq r_0+\sigma k T_0^{\nicefrac{1}{\a}}.$
\end{dem}

\begin{cor} \label{coro}
For every $0<\sigma < 2$ and $\a \in (\nicefrac{1}{2},1)$, let $T_0$ be as in Lemma \ref{it}. Then, for every $u_0$ with compact support, $0\leq u_0 \leq 1$, with $u_0 \neq 0$, there exist $\eps \in (0,1)$ and $b>0$ such that 
$$u(x,t) \geq \eps, \mbox{ for }  T_0\leq t < \tau_{\a} \mbox{ and } \abs x \leq b+\sigma t^{\nicefrac{1}{\a}}.$$
\end{cor}

\begin{dem}
Let $\sigma \in (0,2).$
We have $u(\cdot,\tau)>0$ in $\R^d$ for all $\tau \in [\nicefrac{T_0}{2},\nicefrac{3T_0}{2}]$. Thus, there exists $r_0 >2 d > \tild{C} T_0^{\nicefrac{1}{2\a}}$, where $d$ is a positive constant depending on $T_0$ chosen later, and $\tild{\eps} \in (0,1)$ such that for all $\eps \in (0,\tild{\eps})$:
$$ u(\cdot,\tau)\geq {\eps} \mathds{1}_{B_{r_0}(0)}, \mbox{ in } \R^d$$
Let us define $\underline{u_0}$ by: $\underline{u_0}(y)={\eps} \mathds{1}_{B_{r_0}(0)}(y)$. Thus, $u(\cdot,\tau+t) \geq v(\cdot,t)$, for $t>0$, where $v$ is the solution to \eqref{systeme} with initial condition $\underline{u_0}$ at time $0$. Next, we apply Lemma \ref{it} to $v$:  for all $k\in \N$ such that $kT_0 < \tau_{\a}$:
$$v(x,kT_0) \geq \eps \mbox{  for  } \abs x \leq r_0 + k\sigma T_0^{\nicefrac{1}{\a}} .$$
Consequently:
For all $\tau \in  [\nicefrac{T_0}{2},\nicefrac{3T_0}{2}]$ , for all $k\in \N$ such that $\tau + kT_0< \tau_{\a}$:
$$u(x,\tau +kT_0) \geq \eps \mbox{  for  } \abs x \leq r_0 + k\sigma T_0^{\nicefrac{1}{\a}} .$$
Moreover, the set $\{\tau +kT_0, k \in \N , \tau \in  [\nicefrac{T_0}{2},\nicefrac{3T_0}{2}] \ | \  \tau+ kT_0 < \tau_{\a}\}$ covers all $[\nicefrac{T_0}{2}, \tau_{\a})$.
Let $ T_0 \leq t < \tau_{_a} $, then there exist $\tau  \in  [\nicefrac{T_0}{2},\nicefrac{3T_0}{2}]$ and $ k \in \N $ such that: $t=\tau + kT_0$.\\
Let us define, for $\a$ close to $1$ enough, a constant $d > 0$ independent of $\a$ such that: $(\nicefrac{3T_0}{2} + kT_0)^{\nicefrac{1}{\a}}\leq kT_0^{\nicefrac{1}{\a}}+d$. 
Then:
 $$u(x,t)\geq \eps \mbox{ for } \abs x \leq r_0- \sigma d +\sigma t^{\nicefrac{1}{\a}}.$$
This last statement proves the corollary taking $b=  r_0- 2 d >0$.
\end{dem}

Now, we can prove the first part of Theorem \ref{thm}:
\begin{theo}\label{theo1}
Consider $u$ the solution to \eqref{systeme}, with compactly supported initial datum and $0\leq u_0 \leq 1$, $u_0 \neq 0$. Then:
\begin{itemize}
\item if $\sigma > 2$, then $u(x,t) \rightarrow 0$ uniformly in $\{ \abs x \geq \sigma t^{\nicefrac{1}{\a}}\}$ as $\a \rightarrow 1, t \rightarrow +\infty, t < \tau_{\a}$.
\item if $0< \sigma < 2$, then $u(x,t) \rightarrow 1$ uniformly in $\{ \abs x \leq \sigma t^{\nicefrac{1}{\a}}\}$ as $\a \rightarrow 1, t \rightarrow +\infty, t<\tau_{\a}$.
\end{itemize}
\end{theo}
\begin{dem}
We prove the first statement of the theorem.\\
Let $\sigma$ be such that $\sigma >2$. Recall that $u(x,t) \leq \overline{u}(x,t)$, where:
$$\overline{u}(x,t)=e^t \int_{\R^d} u_0(y) p(x-y,t) dy.$$
Define $R>0$ such that: $\mbox{supp } u_0 \subset B_R(0)$. Let $\a_1 \in (\frac{1}{2},1)$ and $t_0>0$ chosen such that if $x \in \{ x \in \R^d  \ | \ \abs x \geq \sigma t^{\nicefrac{1}{\a}}, \a \in (\a_1,1),t_0\leq  t \leq \tau_{_a}\}$, then $\abs x > 2 R$. Let $x$ be in this set and note that if $\abs{y} \leq R$ then $\abs{y-x} \geq \dis{\frac{\abs x}{2}} \geq \dis{\frac{\sigma t^{\nicefrac{1}{\a}}}{2}} > R$. Moreover, let $\sigma' \in ( 2, \sigma)$, so that: $\abs x -R \geq \sigma ' t^{\nicefrac{1}{\a}}$, for $t$ large enough.Using the fact: $e^{\tau_{\a}}=e^{\frac{\xi_{\a}^{2\a}}{4}}\underset{\a \rightarrow 1}{\sim}C\sin(\a \pi) ^{-1}$,  for $t < \tau_{\a}$, we get: (Recall that the main term in the heat kernel is the Gaussian one)
\begin{eqnarray}
u(x,t) &\leq & Ce^t \left( \int_{\R^d} \frac{t^2 (1-\a)u_0(y)}{\abs{x-y}^{d+4\a}}dy +  \int_{\R^d} \frac{ t \sin (\a \pi)u_0(y)}{\abs{x-y}^{d+2\a}} dy+\int_{\R^d} \frac{u_0(y)e^{-\frac{\abs{y-x}}{4t}^{2\a}}}{t^{\nicefrac{d}{2}}\abs{x-y}^{d(1-\a)}} dy\right)  \nonumber \\
&\leq &  \frac{R^dC(1-\a)}{ \sin(\a \pi) \sigma^{d+4\a}t^{\nicefrac{d}{\a}+2}} + \frac{R^dC}{\sigma^{d+2\a}t^{\nicefrac{d}{\a}+1}}  
+ Ce^t \int_{\abs z \geq \frac{\abs x -R}{(4t)^{\nicefrac{1}{2\a}}}} e^{- \abs z^{2\a}}\abs z^{d(\a-1)}dz \nonumber  \\
&\leq & C t^{-\nicefrac{d}{\a}-2} + C t^{-\nicefrac{d}{\a}-1}+Ce^t \int_{\abs z \geq \frac{\sigma ' t^{\nicefrac{1}{2\a}}}{4^{\nicefrac{1}{2\a}}}} e^{- \abs z^{2\a}}\abs z^{d(\a-1)}dz \nonumber \\
&\leq & C t^{-\nicefrac{d}{\a}-2} + C t^{-\nicefrac{d}{\a}-1}+C e^{t-\frac{ \sigma '^{2\a}}{4} t}\nonumber
\end{eqnarray}
Consequently, since $\sigma> \sigma ' >2$, there exists $\a_2 \in (\a_1,1)$ such that $\sigma ' > 2^{\nicefrac{1}{\a_2}}$. Thus: $(1-\frac{\sigma '^{2\a}}{4})<0, \forall \sigma \in (\a_2,1)$.
Finally, we have:
$$ u(x,t) \rightarrow 0 \mbox{ uniformly in }\{ \abs x \geq \sigma t^{\nicefrac{1}{\a}}\} \mbox{ as }\a \rightarrow 1, t \rightarrow +\infty, t<\tau_{\a}.$$
Now we prove the second statement of the theorem:\\
Given $0<\sigma < 2$, take $\sigma' \in (\sigma,2)$, and apply Corollary \ref{coro} with $\sigma$ replaced by $\sigma'$. Thus, we obtain:
$$-u \leq -\eps \mbox{ in } \omega:=\left\{ (x,t) \in \R^d \times \R^+ \mbox{ such that } T_0 \leq t< \tau_{\a} , \ \abs x \leq b+\sigma' t^{\nicefrac{1}{\a}} \right\},$$
for some $b>0$.
Moreover: $(\partial_t + (-\Delta)^{\a})(1-u)=-u(1-u) \leq - \eps (1-u) \mbox{ in } \omega.$
Let $v$ be the solution to:
\begin{equation} 
\left\{
\begin{array}{rclc}
v_t +(-\Delta)^{\a}v&=&- \eps v,& \quad \R^d, t>0\\
v(y,T_0)&=&1+\frac{e^{\gamma \abs y^{\a}}}{D}  \mathds{1}_{\left\{ \abs y\leq c {\tau_{\a}}^{\nicefrac{1}{\a}}\right\}}(y), &\: \R^d,\\
\end{array}
\right.
\end{equation}
where $\gamma$ and $D$ are constants (depending on $\a$) chosen later,  $c$ is a constant  independent of $\a$ and bigger than $\sigma'$.
There holds:
$$v(x,t)=e^{-\eps (t-T_0)} \left(1+ \dis{\int_{ \abs y \leq c {{\tau_{\a}}}^{\nicefrac{1}{\a}}}\frac{ e^{\gamma \abs y^{\a}}}{D}p(x-y,t-T_0) dy} \right),$$
and $1-u\leq v  \mbox{ in } \omega$.
Now we want to apply lemma 2.1 of \cite{JMRXC}, to have:
$$0\leq 1-u \leq v \mbox{ in } \R^d \times [T_0,\tau_{\a}).$$
Let $w:= 1-u-v$ with initial time $T_0$ and $\abs x \leq r(t):= b+\sigma' t^{\nicefrac{1}{\a}}$.
Let us verify the assumptions of the lemma:
\begin{itemize}
\item Initial datum: $w(.,T_0) \leq 0$ since $1-u \leq 1 \leq v \mbox{ for } t=T_0$
\item Condition outside $\omega$: let $\tau_{\a} > t \geq T_0$ and $\abs x \geq r(t)$ so that $\abs x \leq r(\tau_{\a}).$ We have to verify that $w(x,t) \leq 0$, thus proving that $v(x,t) \geq 1$.
Taking $\a$ closer to $1$ if necessary, we can suppose $r(t) > R$. We use the same inequalities as before:
$$
\begin{disarray}{rcl}
v(x,t) &\geq& C e^{-\eps (t-T_0)}\int_{\tiny{\begin{array}{l}\abs{y-x}\geq \tild{C}  (t-T_0)^{\nicefrac{1}{2\a}}  \\ \abs y \leq c {{\tau_{\a}}}^{\nicefrac{1}{\a}}\end{array}}}\frac{e^{\gamma \abs y^{\a}}}{D} \frac{e^{-\frac{\abs{y-x}}{4 (t-T_0)}^{2\a}}}{ (t-T_0)^{\nicefrac{d}{2}}\abs{y-x}^{d(1-\a)}} dy\\
&\geq& C e^{-\eps (t-T_0)}\int_{\tiny{\begin{array}{l} \abs{z}\geq \nicefrac{\tild{C}}{2^{\nicefrac{1}{\a}}}  \\ \abs {x+ (4 (t-T_0))^{\nicefrac{1}{2\a}} z}\leq c {{\tau_{\a}}}^{\nicefrac{1}{\a}}\end{array}}}\frac{ e^{\gamma   \abs{ x+(4 (t-T_0))^{\nicefrac{1}{2\a}} z}^{\a}}}{D}\frac{e^{-\abs z^{2\a}}}{\abs{z}^{d(1-\a)}} dz\\
&\geq&C e^{-\eps (t-T_0)}e^{\gamma   \abs{ x}^{\a}}\int_{\tiny{\begin{array}{l}\abs{z}\geq \nicefrac{\tild{C}}{4^{\nicefrac{1}{2\a}}}, z.x \geq 0  \\ \abs {x+ (4 (t-T_0))^{\nicefrac{1}{2\a}} z}\leq c {{\tau_{\a}}}^{\nicefrac{1}{\a}}\end{array}}} \frac{e^{-\abs z^{2\a}}}{D\abs{z}^{d(1-\a)}} dz\\
&\geq&C e^{-\eps (t-T_0)}e^{\gamma   (b+ \sigma' t^{\nicefrac{1}{\a}})^{\a}}\int_{\tiny{\begin{array}{l} \abs{z}\geq \nicefrac{\tild{C}}{4^{\nicefrac{1}{2\a}}}, z.x \geq 0  \\ \abs z\leq \frac{c-\sigma '}{2^{\nicefrac{1}{\a}}} {{\tau_{\a}}}^{\nicefrac{1}{2\a}}-\frac{b}{2^{\nicefrac{1}{\a}}}{\tau_{\a}}^{-\nicefrac{1}{2\a}}\end{array}}}\frac{e^{-\abs z^{2\a}}}{D\abs{z}^{d(1-\a)}} dz\\
\end{disarray}
$$
We now choose $D$ so that the right hand side of the above inequality is larger than $e^{(-\eps + \gamma {\sigma'}^{\a})t}$, and we choose $\gamma$ such that: $-\eps + {\sigma'}^{\a} \gamma>0$.

So, we have $v(x,t) \geq 1$ if $T_0$ is large enough, and  as a consequence: $w(x,t) \leq 0$ for  $\tau_{\a} > t \geq T_0$ and $\abs x \geq r(t)$.
\item Let $\tau_{\a} > t \geq T_0$ and $\abs x \leq r(t)$, then we have:
$$w_t(x,t) +(-\Delta)^{\a} w(x,t) \leq - \eps w(x,t),$$
and the last assumption is satisfied.
\end{itemize}
Thus: $w \leq 0$ in $\R \times [T_0, \tau_{\a})$, that is to say:
$$0\leq 1-u(x,t) \leq v(x,t)=e^{-\eps (t-T_0)} \left(1+ \dis{\int_{ \abs y \leq c {{\tau_{\a}}}^{\nicefrac{1}{\a}}}\frac{ e^{\gamma \abs y^{\a}}}{D}p(x-y,t-T_0) dy} \right),$$
for all $ (x,t) \in \R^d \times [T_0, \tau_{\a}).$ 
Finally, we are going to prove that: $v(x,t) \rightarrow 0 \mbox{ uniformly in } \{ \abs x \leq \sigma t^{\nicefrac{1}{\a}}\}$ as $\a \rightarrow 1, t \rightarrow +\infty, t < \tau_{\a}$. For $t  < \tau_{\a}$ and $\abs x \leq \sigma t^{\nicefrac{1}{\a}}$:

$$
\begin{disarray}{rcl}
v(x,t) &\leq &e^{-\eps(t-T_0)}\left(1+\int_{\tiny{\begin{array}{l}\abs y \leq c {{\tau_{\a}}}^{\nicefrac{1}{\a}} \\ \abs{x-y} \leq [ \gamma (t-T_0)]^{\nicefrac{1}{\a}}\end{array}}}e^{\gamma \abs y^{\a}} p(x-y,t-T_0) dy\right.\\
&&\left.\hspace{6cm}+\int_{\tiny{\begin{array}{l}\abs y \leq c {{\tau_{\a}}}^{\nicefrac{1}{\a}} \\ \abs{x-y} \geq [  \gamma (t-T_0)]^{\nicefrac{1}{\a}}\end{array}}}e^{\gamma \abs y^{\a}} p(x-y, t-T_0) dy \right)\\
&\leq & e^{-\eps(t-T_0)}\left(1+\int_{\tiny{\begin{array}{l}\abs y \leq c {{\tau_{\a}}}^{\nicefrac{1}{\a}} \\ \abs{x-y} \leq [  \gamma (t-T_0)]^{\nicefrac{1}{\a}}\end{array}}}\frac{Ce^{\gamma \abs y^{\a}}}{D} \left(\frac{p_{\a}(0)}{ (t-T_0)^{\nicefrac{d}{(2\a)}}}+  \frac{\abs{x-y}}{(t-T_0)^{\nicefrac{(d+1)}{2\a}}}\right) dy\right)\\
&&+ C e^{-\eps(t-T_0)} \left( \int_{\tiny{\begin{array}{l}\abs y \leq c {{\tau_{\a}}}^{\nicefrac{1}{\a}} \\ \abs{x-y} \geq [  \gamma (t-T_0)]^{\nicefrac{1}{\a}}\end{array}}}  e^{\gamma \abs y^{\a}} \frac{(t-T_0)^2(1-\a)}{D \abs{x-y}^{d+4\a}}\right.+ e^{\gamma \abs y^{\a}} \frac{(t-T_0)\sin(\a\pi)}{D\abs{x-y}^{d+2\a}}dy  \\
&&+\left.\int_{\tiny{\begin{array}{l}\abs y \leq c {{\tau_{\a}}}^{\nicefrac{1}{\a}} \\ \abs{x-y} \geq [  \gamma (t-T_0)]^{\nicefrac{1}{\a}}\end{array}}}   \frac{e^{\gamma \abs y^{\a}-\frac{\abs{x-y}^{2\a}}{4 (t-T_0)}}}{ (t-T_0) ^{\nicefrac {d}{2}}\abs{x-y}^{d(1-\a)}}dy\right)\\
&\leq & e^{-\eps(t-T_0)}\left(1+C \left| B_{[\gamma (t-T_0)]^{\nicefrac{1}{\a}}}\right| e^{\gamma \abs x ^{\a}}e^{\gamma^2(t-T_0)} (t-T_0)^{-\nicefrac{(d-1)}{2\a}}+  Ce^{\gamma c^{\a}{\tau_{\a}}} \frac{(1-\a) } {(t-T_0)^2}\right. \\
&&\left.+ Ce^{\gamma c^{\a}{\tau_{\a}}} \frac{\sin(\a\pi)}{ (t- T_0)}+\int_{\tiny{\begin{array}{l}\abs {x+(4(t-T_0))^{\nicefrac{1}{2\a}}z} \leq c {{\tau_{\a}}}^{\nicefrac{1}{\a}} \\ \abs{z} \geq [ \frac{\gamma}{2} ]^{\nicefrac{1}{\a}}(t-T_0)^{\nicefrac{1}{2\a}}\geq 1\end{array}}}  e^{\gamma \abs x^{\a}}e^{2 \gamma \sqrt{(t-T_0)}\abs z^{\a}-\abs{z}^{2\a}}dz\right)\\
&\leq& e^{-\eps(t-T_0)}+ C (t-T_0)^{\nicefrac{(d+1)}{2\a}} e^{-\eps(t-T_0)} e^{\gamma \sigma^{\a}t}e^{\gamma^2(t-T_0)} + Ce^{-\eps(t-T_0)}e^{(-1+\gamma c^{\a}){\tau_{\a}}} \\
&&+e^{-\eps(t-T_0)}e^{\gamma \sigma^{\a}t}\int_{\R^d}e^{-(\abs z^{\a} - \gamma \sqrt{(t-T_0)})^2}e^{\gamma^2(t-T_0)}dz\\
&\leq&e^{-\eps(t-T_0)}+ C e^{(-\eps+\sigma^{\a}\gamma +\gamma^2)(t-T_0)}    (t-T_0)^{\nicefrac{(d+1)}{2\a}}e^{\gamma \sigma^{\a}T_0}  + Ce^{-\eps(t-T_0)} e^{(-1+\gamma c^{\a}){\tau_{\a}}}\\
&&+e^{(-\eps+\sigma^{\a}\gamma +\gamma^2)(t-T_0)} e^{\gamma \sigma^{\a}T_0}\int_{\R^d}  e^{-\abs z^{2\a}}dz\\
\end{disarray}
$$
We have the result if $\gamma$ is chosen so that: $-\eps+\sigma^{\a}\gamma +\gamma^2<0 $, that is to say: $\gamma < \gamma_1= \frac{-\sigma^{\a}+\sqrt{\sigma^{2\a}+4 \eps}}{2} \underset{\eps \rightarrow 0}{\sim} \frac{\eps}{\sigma^{\a}}.$ With such a $\gamma$, we have for $\eps$ small enough: $-1+c^{\a}\gamma<0$.
Eventually, taking $\eps$ smaller if necessary, we have $\eps \sigma'^{-\a} < \gamma_1$, and consequently, we take: $ \eps \sigma'^{-\a} < \gamma < \gamma_1$ and we get:
$$u(x,t) \rightarrow 1 \mbox{ uniformly in } \{ \abs x \leq \sigma t^{\nicefrac{1}{\a}}\} \mbox{ as } \a \rightarrow 1, t \rightarrow +\infty, t<\tau_{\a}.$$
\end{dem}

\subsection{The exponential propagation phase}\label{propexp1}
Let us now worry about the behaviour of \eqref{systeme} for $t \geq \tau_{\a}$. The initial condition is $u(x, \tau_{\a})= \int_{\R^d}u_0(y) p(x-y,\tau_{\a})dy$.
To prove the second part of Theorem \ref{thm}, we use the argument developped in \cite{JMRXC}. However it is not sufficient to give us the evolution of the level set $\{ x \in \R^d \  | \ u(x,t)=\eps \}$ with $\eps>0$ independent of $\a$: indeed the fundamental solution contains a $\sin(\a \pi)$ which, a priori, only tells us something about values of $u$ of the magnitude $\eps_{\a}=\sin(\a \pi)^{1+ \kappa}, \kappa >0$. So the proof will follow the iterative scheme of \cite{JMRXC} but, within each iteration, two steps will be needed: the first one will study the evolution of the level set  $\{ x \in \R^d \  | \ u(x,t)=\eps_{\a} \}$ and the second one will use the intermediate Proposition \ref{youpi}, to pass from $\eps_{\a}$ to $\eps$ (independent of $\a$).
\begin{lemme}\label{ite2}
For every $0<\sigma<\frac{1}{d+2\a}$ and $\a \in (\nicefrac{1}{2},1)$, there exist $\eps_{0} \in (0,1)$, $T_{\a} > 1$ depending on $\sigma$ and $\eps_0$ and of order $\tau_{\a}$, and $\tild{\tau_{\a}}>0$ for which the following holds.
Given $r_0 >\tau_{\a}$, $\underline{\eps} \in (0, \eps_{0})$, $\eps_{\a} \in (0, \underline{\eps})$, $a_{0,\a}$ be defined by $a_{0,\a} r_0^{-(d+2\a)}=\eps_{\a}$, and let 
$$
\underline{u_{0,\a}}(x)= \left\{
\begin{array}{ll}
a_{0,\a} r_0^{-(d+2\a)} & \mbox{ if $\abs x \leq r_0$} \\
a_{0,\a} \abs x^{-(d+2\a)} & \mbox{ if $\abs x \geq r_0$} 
\end{array}
\right.
$$
Then, the solution to $u_t+(-\Delta)^{\a}u=u-u^2$ with initial condition $\underline{u_{0,\a}}$ and initial time $\tau_{\a}$ satisfies, for all $k\in \N$:
$$u(x, \tau_{\a}+\tild{\tau_{\a}}+(k+1)T_{\a}) \geq \underline{ \eps} \  \mbox{  for  } \ \abs x \leq (r_0 -M) e ^{\sigma k T_{\a}} ,$$
where $M$ is defined in section \ref{inter}. (Recall that $M$ is large enough so that the principal Dirichlet eigenvalue of $(-\Delta)^{\a}-I$ in $B_M$ is negative.)
\end{lemme}
\begin{dem}
 The difference with the preceding paragraph is that, here, the emphasis is laid on the $\dis{\frac{\sin(\a \pi) t}{\abs x^{d+2\a}}}$ term in the heat kernel. \\
Let $\kappa >0$ be a constant independent of $\a$, large enough and $\a_1 \in (\nicefrac{1}{2},1)$ such that: $\eps_{\a} :=\sin(\a \pi)^{1+\kappa}, \ \forall \a \in (\a_1,1).$ 

 Define $\delta_{\a}= \sqrt{\eps_{\a}}$ and $T_{\a} = \ln\left( \dis{\frac{\delta_{\a}}{B\eps_{\a}}}\right)$ so that (see \ref{T_0}):
\begin{equation}\label{T2}
\forall t\in (0,T_{\a}), \forall x \in \R^d, \quad u(x ,\tau_{\a}+t) \leq \delta_{\a}.
\end{equation}
Notice that $T_{\a}$ is of order $\tau_{\a}$, and thus: $r_0 > \tau_{\a} > \tild{C} T_{\a}^{\nicefrac{1}{2\a}}$, for $\a$ close to $1$ enough, where $\tild{C}$ is defined so that $ \dis{\frac{C(1-\a) t^2}{\abs{y-x}^{d+4\a}}} \leq \dis{\frac{\sin(\a \pi) t}{\abs{y-x}^{d+2\a}}}$, if $\abs{x-y} \geq \tild{C} t^{\nicefrac{1}{2\a}}$. (see \eqref{neg}). \\
Now, we prove the lemma:\\
\textbf{For $k=0$:} We get, for all $x \in \R^d$ and $t\in (0,T_{\a})$:
\begin{eqnarray} \label{dec}
u(x,\tau_{\a}+t) & \geq & e^{(1-\delta_{\a})t}\int_{\abs{y-x}\leq \tild{C} t^{\nicefrac{1}{2\a}}}\hat{C} t^{-\nicefrac{(d+1)}{2\a}}\underline{u_{0,\a}}(y)dy\nonumber\\
&&+C e^{(1-\delta_{\a})t} \int_{\abs{x-y}\geq \tild{C}t^{\nicefrac{1}{2\a}}}  \frac{\sin(\a\pi)t}{ \abs{x-y}^{d+2\a}} \underline{u_{0,\a}}(y) dy\nonumber\\
&:=&w(x, \tau_{\a}+t) 
\end{eqnarray}
(Notice that we have here removed the Gaussian part of the heat kernel.) As in Section \ref{lin}, we get: $w$ is a positive, radially symmetric and nonincreasing function.
And, for $\abs x \leq r_{0}$  :
\begin{eqnarray*}
u(x,\tau_{\a}+T_{\a}) &\geq&w(x, \tau_{\a}+T_{\a}) \nonumber\\
&\geq&w(r_0e_1, \tau_{\a}+T_{\a}) \nonumber\\
& \geq &Ce^{(1-\delta_{\a})T_{\a}}  \int_{\tiny{\begin{array}{l}\abs{r_0e_1-y}\geq \tild{C}T_0^{\nicefrac{1}{2\a}}\\ \abs y \leq r_0 \end{array}}} \frac{\sin(\a\pi)}{\abs{r_0e_1-y}^{d+2\a}}\underline{u_{0,\a}}(y) dy \nonumber\\
& \geq &C\eps_{\a}\frac{ \sin(\a \pi)}{ \eps_{\a}^{\frac{1-\delta_{\a}}{2}}\tild{C}^{d+2\a}T_{\a}^{\nicefrac{d}{2\a}+1}} \nonumber\\
&\geq& C \eps_{\a} \frac{ \sin(\a \pi)^{1-\frac{1-\delta_{\a}}{2}(1+\kappa)}}{ T_{\a}^{\nicefrac{d}{2\a}+1}} \nonumber \\
& \geq & \eps_{\a} \nonumber
\end{eqnarray*}
The last inequality is obtained taking $\kappa$ large enough.
Consequently, $w$ is bigger than the solution $v$ to \eqref{v} with initial condition $\eps_{\a} \mathds{1}_{B_{M-1}[(r_0-M)e_1]} $ at time $\tau_{\a}+T_{\a}$. Note that $r_0>\tau_{\a}>M$ for $\a$ close to $1$ enough.
Then, we use Proposition \ref{youpi} to the solution $v$, and we get the existence of $r \in B_{M}[(r_0-M)e_1]$ such that: $$w(r,\tau_{\a}+T_{\a}+\tild{\tau_{\a}}) \geq v(r,\tild{\tau_{\a}})\geq \underline{\eps}.$$
Finally, since $w$ is radially symmetric and nonincreasing, we get:
$$u(x, \tau_{\a}+T_{\a}+ \tild{\tau_{\a}}) \geq \underline{\eps}, \ \forall \abs x \leq r_0-M.$$
\textbf{For $k=1$:} First, we look for $x_1$ such that: 
$$ u(x,\tau_{\a}+ T_{\a}) \geq \sin(\a \pi)^{1+\kappa}=\eps_{\a}, \ \mbox{ for } \abs x \leq x_1.$$ 
For every $\delta \in (0,1)$, $\delta >> \underline{\eps}$, we have that inequality \eqref{dec} for $t=T_{\a}$ and  $\abs x \geq r_{0}$  leads to:
$$
\begin{disarray}{rcl}
u(x,\tau_{\a}+T_{\a}) &\geq&w(x, \tau_{\a}+T_{\a}) \nonumber\\
& \geq &Ce^{(1-\delta_{\a})T_{\a}} \int_{\tiny{\begin{array}{l}\tild{C}T_{\a}\geq \abs{x-y} \geq \tild{C}T_{\a}^{\nicefrac{1}{2\a}}\\ \abs y \leq \abs x \end{array}}} \frac{\sin(\a\pi)}{ \abs{x-y}^{d+2\a}}\underline{u_{0,\a}}(y) dy \nonumber\\
& \geq & C\frac{a_{0,\a}e^{(1-\delta_{\a}-\delta)T_{\a}}\sin(\a\pi)^{1-\frac{\delta(1+\kappa)}{2}}}{\abs x ^{d+2\a} \tild{C}^{d+2\a}T_{\a}^{d+2\a}} \nonumber\\
& \geq & C\frac{a_{0,\a}e^{(1-\delta_{\a}-\delta)T_{\a}}}{\abs x ^{d+2\a} } \nonumber\\,
\end{disarray}
$$
the last inequality is obtained for $\kappa >\frac{2}{\delta}-1$.
Let  us define $x_1$ by:
 $$  C\frac{a_{0,\a}e^{(1-\delta_{\a}-\delta)T_{\a}}}{ {x_1} ^{d+2\a} } =\eps_{\a}$$
Since $a_{0,\a}= \eps_{\a} r_0^{d+2\a}$, we get:
$$x_1=r_0 C e^{\frac{1-\delta_{\a}-\delta}{d+2\a}T_{\a}}.$$
Consequently, for each $\a<1$, $0<\sigma < \frac{1}{d+2\a}$, we take $ \delta_{\a}$ and $\delta$ small enough so that:  
\begin{equation}\label{delta}
\forall \a \in (\a_1,1), \  \sigma <\frac{1-\delta_{\a}-\delta}{d+2\a}< \frac{1}{d+2\a}.
\end{equation} 
Thus, taking $\a$ closer to $1$ if necessary, we have: $$C e^{\frac{1-\delta_{\a}-\delta}{d+2\a}T_{\a}} \geq e^{\sigma T_{\a}}$$
Now, let us define: $\overline{r_1}= r_0 e^{\sigma T_{\a}}$ so that $\overline{r_1} < x_1$, and  since $w$ is positive, radially symmetric and nonincreasing, we obtain:
$$u(x, \tau_{\a}+T_{\a}) \geq \eps_{\a }, \mbox{ for } \abs x \leq \overline{ r_1}.$$
Furthermore, we have:
$u(x, \tau_{\a}+T_{\a}) \geq \dis{ \frac{a_{1,\a}}{\abs x^{d+2\a}}} \mbox{ for } \abs x \geq \overline{r_1},$ where  $ a_{1,\a} = \eps_{\a} {x_1}^{d+2\a}$.

Thus:
$$u(\cdot,\tau_{\a}+T_{\a}) \geq \underline{u_{1,\a}},$$
where $\underline{u_{1,\a}}$ is given by the same expression as $\underline{u_{0,\a}}$ with $(r_0,a_{0,\a})$ replaced by  $(\overline{r_1},a_{1,\a})$.

Finally, we use the case $k=0$ to make the connection with the level set $\{x \in \R^d \ | \ u(x,t)=\underline{\eps} \}$, replacing the initial time $\tau_{\a}$ by $\tau_{\a}+T_{\a}$ and the initial condition by $\underline{u_{1,\a}}$, to get: 
$$u(x,\tau_{\a}+ \tild{\tau_{\a}}+2T_{\a}) \geq \underline{\eps}, \mbox{ for } \abs x \leq r_1:=(r_0-M)e^{\sigma T_{\a}}.$$
We can repeat the argument above, to get :
$$u(x, \tau_{\a}+ \tild{\tau_{\a}}+(k+1)T_{\a}) \geq \eps_{\a}, \mbox{ for } \abs x \leq r_k,$$
for all $k \in \N$, with $r_k \geq (r_0-M)e^{\sigma kT_{\a}}$.
\end{dem}
\begin{cor} \label{coro2}
For every $0<\sigma < \frac{1}{d+2\a}$ and $\a \in (\nicefrac{1}{2},1)$, let $T_0$ and $\tild{\tau_{\a}}$ be as in Lemma \ref{ite2}. Then, for every $u_0$ with compact support, $0\leq u_0 \leq 1$, with $u_0 \neq 0$, there exist $\underline{\eps} \in (0,1)$, $\overline{C}>0$ a constant independent of $\a$, and $b_{\a}>0$ such that 
$$u(x,t) \geq \underline{\eps}, \mbox{ if }
t \geq \overline{C} \tau_{\a} \mbox{ and } \abs x \leq b_{\a} e^{\sigma t},$$
where $b_{\a}$ is proportional to $e^{-\underline{C}\sigma \tau_{\a}}$, $ \underline{C}$ is a constant independent of $\a$ and strictly smaller than $\overline{C}$.
\end{cor}
\begin{dem}
From Lemma \ref{ite2}, we have $\delta$ defined by \eqref{delta}. Let $\kappa >\frac{2}{\delta}-1$ and $\a_1 \in (\nicefrac{1}{2},1)$ be such that: $\eps_{\a} :=\sin(\a \pi)^{1+\kappa}, \ \forall \a \in (\a_1,1)$.\\
Using the proof done for $t<\tau_{\a}$, we get the existence of $\eps \in (0,1)$ so that, for $\a$ closer to $1$ if necessary:
$$u(x, \tau_{\a}) \geq \eps \geq a_{0,\a} , \mbox{ for } \abs x \leq  r_{\a},$$
where $r_{\a} > \tau_{\a}^{\nicefrac{1}{\a}}$, and $a_{0,\a}= \eps_{\a}  r_{\a}^{d+2\a}$. 
Thus $$u(\cdot, \tau_{\a}) \geq a_{0,\a} \mathds{1}_{B_{r_{\a}}(0)} \mbox{ in } \R^d,$$
and $u(\cdot,  \tau_{\a}+t) \geq v(\cdot,t)$ for t>0, where $v$ is the solution to \eqref{systeme} with initial condition $a_{0,\a} \mathds{1}_{B_{r_{\a}}(0)}$. 
Next, we denote by $\overline{T_0}$ the time before which the solution $u$ reaches $\delta$.
Inequality \eqref{T_0} leads to  $\overline{T_0}= \ln(\frac{\delta}{B \eps_{\a}})$,consequently:
$$\forall t \in ( \tau_{\a},  \tau_{\a}+\overline{T_0}), \ \ u(\cdot,t) \leq \delta.$$
Note that $\overline{T_0 }$ is of order $\tau_{\a}$ and for $\a$ close to $1$ enough: $T_{\a} < \frac{2}{3} \overline{T_0}$.
Since $u\leq \delta$, we get, for $ t \in [\nicefrac{\overline{T_0}}{3},\nicefrac{\overline{T_0}}{3}+T_{\a}]$:
$$
\begin{disarray}{rcl}
u( x,  \tau_{\a}+t) &\geq & e^{(1-\delta)t}\int_{\R^d} p(x-y,t) a_{0,\a} \mathds{1}_{B_{r_{\a}}(0)}(y) dy \\
&\geq& C \int_{\tiny{ \begin{array}{l} \abs{x-y} \geq \tild{C} t^{\nicefrac{1}{2\a}} \\ \abs y \leq r_{\a} \end{array}}}e^{\frac{1-\delta}{3}\overline{T_0}} \frac{a_{0,\a} t \sin(\a \pi) } {\abs{x-y}^{d+2\a}}dy\\
&\geq&C\int_{\tiny{ \begin{array}{l} \abs{x-y} \geq \tild{C}t^{\nicefrac{1}{2\a}} \\ \abs y \leq 1 \end{array}}}\left(\frac{\delta}{B \eps_{\a}}\right)^{\frac{1-\delta}{3}}\frac{a_{0,\a} \sin(\a \pi) } {\abs{x-y}^{d+2\a}}dy:= \tild{w}(x,\tau_{\a}+ t)\\
\end{disarray}
$$
Moreover, using the fact that for $\abs x \geq r_{\a} > \tau_{\a}^{\nicefrac{1}{\a}}$ and $\abs y \leq r_{\a}$: $\abs{x-y} \leq 2 \abs x$, we obtain:
$$
\begin{disarray}{rcl}
u( x,  \tau_{\a}+t)&\geq&C \int_{\abs y \leq 1} \frac{ a_{0,\a} } {\sin(\a \pi)^{(1+\kappa)\frac{1-\delta}{3}-1}\abs{x}^{d+2\a}}dy\\
&\geq& C\frac{ a_{0,\a} } {\sin(\a \pi)^{\kappa_0}\abs{x}^{d+2\a}},
\end{disarray}
$$
where $\kappa_0  >0$.
Taking $\a$ closer to $1$ if necessary, we get, for $ t \in [\nicefrac{\overline{T_0}}{3},\nicefrac{\overline{T_0}}{3}+T_{\a}]$ and $x \in \R^d$:
$$u( x,  \tau_{\a}+t) \geq\dis{ \frac{a_{0,\a}}{\abs x^{d+2\a}}}, \mbox{ for } \abs x \geq r_{\a}.$$
As a consequence, using the fact $\tilde{w}$ is symmetric radially nonincreasing:
 $$u( x,   \tau_{\a}+t) \geq \tild{w}(x, \tau_{\a} + t) \geq \tild{w}(r_{\a}e_1 , \tau_{\a} + t)\geq \eps_{\a}, \mbox{ for } \abs x \leq r_{\a}.$$
Finally:
$$ u(\cdot , \tau_{\a}+t) \geq \underline{ u_{0,\a}}, \ \ \forall  t \in [\nicefrac{\overline{T_0}}{3},\nicefrac{\overline{T_0}}{3}+T_{\a}],$$
where $\underline{u_{0,\a}}$ is the initial condition in Lemma \ref{ite2}, with $r_0$ replaced by $r_{\a}$.

Next, we can apply Lemma \ref{ite2} to the solution $ u(\cdot ,\cdot+\tau_0)$ for all  
$\tau_0 \in [ \nicefrac{\overline{T_0}}{3},\nicefrac{\overline{T_0}}{3}+T_{\a}]$.
Indeed, $\{ \tau_{\a}+\tau_0+(k+1)T_{\a}+\tild{\tau_{\a}}, \ k \in \N , \tau_0 \in     [ \nicefrac{\overline{T_0}}{3},\nicefrac{\overline{T_0}}{3}+T_{\a}]  \}$ covers all $(\tild{\tau_{\a}}+\tau_{\a}+ \nicefrac{\overline{T_0}}{3}+T_{\a}, + \infty)$.
Let $\overline{C} $ be a constant such that  $\tau_{\a}+ \nicefrac{\overline{T_0}}{3}+T_{\a}+\tild {\tau_{\a}} \leq \overline{C}\tau_{\a}$. If $t \geq \overline{C}\tau_{\a}$, then there exist $\tau_0 \in   [ \nicefrac{\overline{T_0}}{3},\nicefrac{\overline{T_0}}{3}+T_{\a}]$ and $k \in \N$ such that: $$t= \tau_{\a}+\tau_0+(k+1)T_{\a}+\tild{\tau_{\a}}.$$

Then: $$u(x,t) \geq \underline{\eps}, \mbox{ if }  t \geq  \overline{C}\tau_{\a} \mbox{ and } \abs x \leq b_{\a} e^{\sigma t},$$
with : $b_{\a}=C \ e^{-\sigma \underline{C} \tau_{\a}} > 0$, where $C>0$ is a constant independent of $\a$ and $\underline{C}$ is a constant independent of $\a$ such that $\tau_{\a}+ \nicefrac{\overline{T_0}}{3}+2T_{\a} +\tild{\tau_{\a}}\leq \underline{C} \tau_{\a}$. Taking $\overline{C}$ larger if necessary, we can assume:    $\underline{C} < \overline{C}$.
\end{dem}

Now, we can prove the second part of Theorem \ref{thm}:
\begin{theo}
Under the assumptions of Theorem \ref{theo1}, there exists a constant $\overline{C} >0$ such that:
\begin{itemize}
\item if $\sigma > \dis{\frac{1}{d+2\a}}$, then $u(x,t) \rightarrow 0$ uniformly in $\{ \abs x \geq e^{\sigma t}\}$ as $\a \rightarrow 1, t \rightarrow +\infty, t > \overline{C} \tau_{\a}$.
\item if $0< \sigma < \dis{\frac{1}{d+2\a}}$, then $u(x,t) \rightarrow 1$ uniformly in $\{ \abs x \leq e^{\sigma t}\}$ as $\a \rightarrow 1, t \rightarrow +\infty, t >\overline{C}  \tau_{\a}$.
\end{itemize}
\end{theo}

\begin{dem}
We prove the first statement of the theorem.
Let $\sigma$ be such that $\sigma >\dis{\frac{1}{d+2\a}}$, and $x$ such that $ \abs x \geq e^{\sigma t}$.
Recall that $u(x,t) \leq \overline{u}(x,t)$, where:
$$\overline{u}(x,t)=e^t \int_{\R^d} u_0(y) p(x-y,t) dy.$$
Define $R>0$ such that: $\mbox{ supp } u_0 \subset B_R(0)$. For, $ \abs y \leq R$,  we take $t$ large enough to have $\abs{ x-y} \geq \dis{\frac{e^{\sigma t }}{2}}$. Then:
$$
\begin{disarray}{rcl}
u(x,t)&\leq&Ce^t \left( \int_{\abs y \leq R} \frac{(1-\a)t^2}{\abs{x-y}^{d+4\a}}dy +\int_{\abs y \leq R} \frac{\sin(\a \pi)  t}{\abs{x-y}^{d+2\a}}dy +\int_{\abs y \leq R} \frac{e^{-\frac{\abs{x-y}}{4t}^{2\a}}}{t^{\nicefrac{d}{2}}\abs{x-y}^{d(1-\a)}}dy\right)\\
&\leq&C(1-\a)R^d t^2 e^{(1-\sigma(1+4\a))t}+C R^d t e^{(1-\sigma(1+2\a))t}+CR^de^{(1-(1-\a)\sigma)t-\frac{e^{2\sigma \a t}}{4t}} t^{\nicefrac{d} {2}}\\
\end{disarray}
$$
Hence, for $\sigma >\dis{\frac{1}{d+2\a}}$, we obtain:
$$ u(x,t) \rightarrow 0 \mbox{ uniformly in } \{ \abs x \geq e^{\sigma t} \mbox{ as  } \a \rightarrow 1, t \rightarrow +\infty, t > \tau_{\a}\}.$$
Now, we prove the second statement of the theorem. The argument parallels that of Theorem \ref{theo1}. Using the proof done for $t < \tau_{\a}$, we know there exists $\eps \in (0,1)$ such that:$$u(\cdot, \tau_{\a}) \geq \eps \mathds{1}_{B_{r_{\a}}(0)},$$
where $r_{\a}$ is smaller than $2 \tau_{\a}^{\nicefrac{1}{\a}}$. 
As in the beginning of the proof of  Corollary \ref{coro2}, we have, for $ t \in  [\nicefrac{\overline{T_0}}{3},\nicefrac{\overline{T_0}}{3}+T_{\a}]$:
$$u(x,t) \geq u_{0,\a} = \left\{ \begin{array}{rl}\frac{a_{0,\a}}{\abs x ^{d+2\a}}, & \abs x \geq r_{\a}\\
\eps_{\a}, & \abs x \leq r_{\a} \end{array} \right. . $$
Given $0<\sigma <\dis{\frac{1}{d+2\a}}$, take $\sigma'$ with: $0<\sigma < \sigma '<\dis{\frac{1}{d+2\a}}$ and  apply Corollary \ref{coro2} with $\sigma$ replaced by $\sigma'$.
Thus, we obtain:
$$-u \leq -\underline{\eps} \mbox{ in } \omega :=\left\{(x,t) \in \R^d \times \R^+ \quad  | \quad t\geq \overline{C}\tau_{\a}, \abs x \leq b_{\a}e^{\sigma' t} \right\},$$
where  $b_{\a}=(r_{\a}-M) \ e^{-\sigma' \underline{C} \tau_{\a}} > 0$.
What's more: $(\partial_t + (-\Delta)^{\a})(1-u)=-u(1-u) \leq - \eps (1-u) \mbox{ in } \omega.$
Define $v$ the solution to:
\begin{equation} 
\left\{
\begin{array}{rclc}
v_t +(-\Delta)^{\a}v&=&- \underline{\eps} v,& \quad \R^d, t> \overline{C}\tau_{\a}\\
v(y, \overline{C}\tau_{\a})&=&1+\dis{\frac{e^{-\gamma\sigma\overline{C}\tau_{\a}}\abs y^{\gamma}}{D}},& \: \R^d\\
\end{array}
\right.
\end{equation}
where $\gamma \in (0,2\a)$  and $D$ are constants, independent of $\a$, chosen later.
This solution is given by: $$v(x,t)=e^{-\eps (t- \overline{C}\tau_{\a})}\left(1+ \dis{\int_{ y \in \R^d}\frac{e^{-\gamma\sigma\overline{C}\tau_{\a}} \abs y^{\gamma}}{D }p(x-y,t- \overline{C}\tau_{\a}) dy} \right),$$
and $1-u\leq v  \mbox{ in } \omega$. To get:
$$0\leq 1-u \leq v \mbox{ in } \R \times (\overline{C}\tau_{\a},+\infty),$$
let us verify the assumptions of the lemma 2.1 in \cite{JMRXC}. Let $w:= 1-u-v$ with initial time $\overline{C}\tau_{\a}$ and $\abs x \leq r(t):= b_{\a}e^{\sigma' t}$. Remember that in this case: $b_{\a}=C \ e^{-\sigma' \underline{C} \tau_{\a}} > 0$.
\begin{itemize}
\item Initial datum: $w(.,\overline{C}\tau_{\a}) \leq 0$ since $1-u \leq 1 \leq v \mbox{ for } t=\overline{C}\tau_{\a}$
\item Condition outside $\omega$: let $ t \geq \overline{C}\tau_{\a}$ and $\abs x \geq r(t)$ . We have to verify that $w(x,t) \leq 0$, proving that $v(x,t) \geq 1$.
Taking  $\a$ closer to $1$ if necessary, we can suppose $r(t) > R$. We use the same inequalities as before  taking $\overline{C}$ larger if necessary and using the fact that $\sigma < \sigma '$:
$$
\begin{disarray}{rcl}
v(x,t)&\geq& e^{-\underline{\eps}(t-\overline{C}\tau_{\a})}\int_{y \in \R^d}  \frac{ e^{-\gamma \sigma \overline{C}\tau_{\a}}\abs y^{\gamma}}{D} p(x-y, t-\overline{C}\tau_{\a})dy\\
&\geq &C e^{-\underline{\eps}(t-\overline{C}\tau_{\a})}\int_{\abs{x-y}\geq \tild{C}(t-\overline{C}\tau_{\a})^{\nicefrac{1}{2\a}}}\frac{\sin(\a\pi)(t-\overline{C}\tau_{\a})   e^{-\gamma \sigma \overline{C}\tau_{\a}} \abs y^{\gamma}}{D  \abs{x-y}^{d+2\a}} dy\\
&\geq &C e^{-\underline{\eps}(t-\overline{C}\tau_{\a})}\int_{\nicefrac{\abs x }{2}\geq\abs{x-y}\geq \tild{C}(t-\overline{C}\tau_{\a})^{\nicefrac{1}{2\a}}}\frac{\sin(\a\pi) (t-\overline{C}\tau_{\a})   e^{-\gamma \sigma \overline{C}\tau_{\a}}\abs y^{\gamma}}{D  \abs{x-y}^{d+2\a}} dy\\
&\geq& C e^{-\underline{\eps}(t-\overline{C}\tau_{\a})}\int_{ (t-\overline{C}\tau_{\a})^{-\nicefrac{1}{2\a}}\nicefrac{\abs x}{2} \geq \abs{z}\geq \tild{C}}\frac{   e^{-\gamma \sigma \overline{C}\tau_{\a}}\abs x^{\gamma}\sin(\a\pi) }{D\abs{z}^{d+2\a}} dz\\
&\geq& CD^{-1}e^{-\underline{\eps}(t-\overline{C}\tau_{\a})}b_{\a}^{\gamma} e^{-\gamma \sigma \overline{C}\tau_{\a}}e^{\sigma ' \gamma t}\sin(\a\pi)\\
&\geq & e^{(-\underline{\eps} +\gamma \sigma') ( t-\overline{C}\tau_{\a})} e^{\gamma \sigma ' \overline{C}\tau_{\a}-\gamma \sigma \overline{C}\tau_{\a}-\gamma  \sigma ' \underline{C}\tau_{\a}-\tau_{\a}}\\
&\geq & e^{(-\underline{\eps} +\gamma \sigma')(t-\overline{C}\tau_{\a})},
\end{disarray}
$$
these inequalities are obtained taking $0<D\leq C$, independent of $\a$.
Thus, if $\gamma$ is chosen so that: $-\underline{\eps} +\gamma \sigma' > 0$:
$$v(x,t) \geq 1 \geq 1-u(x,t) , \mbox{ for } t \geq \overline{C} \tau_{\a} \mbox{ and } \abs x \geq r(t).$$
Notice once again that we have removed the Gaussian term.
\item Let $ t \geq \overline{C}\tau_{\a}$ and $\abs x \leq r(t)$, then we have:
$$w_t(x,t) +(-\Delta)^{\a} w(x,t) \leq - \underline{\eps} w(x,t),$$
and the last assumption is satisfied.
\end{itemize}
Now, by Lemma 2.1 in \cite{JMRXC}, we have: $w \leq 0$ in $\R^d \times [\overline{C}\tau_{\a}, +\infty)$, that is to say:
$$0\leq 1-u(x,t) \leq v(x,t)=e^{-\eps (t-\overline{C}\tau_{\a})} \left(1+ \dis{\int_{  y \in \R^d }}\frac{e^{-\gamma \sigma \overline{C}\tau_{\a}} \abs y^{\gamma}}{D}p(x-y,t-\overline{C}\tau_{\a}) dy \right),$$
for all $ (x,t) \in \R^d \times [\overline{C}\tau_{\a},+\infty).$
Finally, we are going to prove that: $v(x,t) \rightarrow 0$ uniformly in $\{ \abs x \leq e^{\sigma t}\}$ as $\a \rightarrow 1, t \rightarrow +\infty, t >\overline{C}\tau_{\a}$:

$$
\begin{disarray}{rcl}
v(x,t) &\leq& Ce^{-\underline{\eps} (t-\overline{C}\tau_{\a})} \left(1+\int_{\abs{x-y}\leq 1} \frac{e^{-\gamma \sigma \overline{C}\tau_{\a}} \abs y^{\gamma}}{D }dy \right.+\int_{\abs{x-y} \geq 1} \frac{(1-\a)(t-\overline{C}\tau_{\a})^2e^{-\gamma \sigma \overline{C}\tau_{\a}} \abs y^{\gamma}}{D \abs{x-y}^{d+4\a}}dy \\
&&+\int_{\abs{x-y} \geq 1} \frac{\sin(\a\pi) (t-\overline{C}\tau_{\a}) e^{-\gamma \sigma \overline{C}\tau_{\a}} \abs y^{\gamma}}{D   \abs{x-y}^{d+2\a}}dy \left.+\int_{\abs{x-y} \geq 1} \frac{e^{-\frac{\abs{x-y}^{2\a}}{4(t-\overline{C}\tau_{\a})}} e^{-\gamma \sigma \overline{C}\tau_{\a}} \abs y^{\gamma}}{D  (t-\overline{C}\tau_{\a})^{\nicefrac{d}{2}} \abs{x-y}^{d(1-\a)}}dy\right)\\
&\leq& C e^{-\underline{\eps} (t-\overline{C}\tau_{\a})} \left(1+ \frac{e^{-\gamma \sigma \overline{C}\tau_{\a}} (\abs x^{\gamma}+1)}{D}+ \int_{\abs z \geq 1} \frac{(1-\a) e^{-\gamma \sigma \overline{C}\tau_{\a}}(t-\overline{C}\tau_{\a})^2 (\abs x^{\gamma}+ \abs z ^{\gamma})}{D \abs{z}^{d+4\a}} dz \right.\\
&&+  \left. \int_{\abs{z} \geq 1} \left(\frac{\sin(\a \pi)(t-\overline{C}\tau_{\a})}{D   \abs{z}^{d+2\a}}+ \frac{e^{-\frac{\abs{z}^{2\a}}{4(t-\overline{C}\tau_{\a})}}}{D \abs{z}^{d(1-\a)}}\right)  e^{-\gamma \sigma \overline{C}\tau_{\a}} (\abs x^{\gamma}+ \abs z ^{\gamma})dz\right)\\
&\leq& C e^{-\underline{\eps} (t-\overline{C}\tau_{\a})}\left(1+\int_{\abs z \geq 1}\frac{(t-\overline{C}\tau_{\a})^2}{ \abs z^{ d + 4\a -\gamma}}+\frac{(t-\overline{C}\tau_{\a})}{\abs z^{  d + 2\a -\gamma}}dz+ \int_{\R^d} e^{-\frac{\abs{z}^{2\a}}{4(t-\overline{C}\tau_{\a})}}dz\right)\\
&&+Ce^{(-\underline{\eps}+\gamma \sigma)(t-\overline{C}\tau_{\a})} \left(1+\int_{\abs z \geq 1} \frac{(t-\overline{C}\tau_{\a})^2}{ \abs z^{ d + 4\a }} dz+ \int_{\abs z \geq 1}\frac{(t-\overline{C}\tau_{\a})}{ \abs z^{ d + 2\a}}dz  +\int_{\R^d}e^{-\frac{\abs{z}^{2\a}}{4(t-\overline{C}\tau_{\a})}}dz\right).
\end{disarray}
$$
Notice that all the integrals converge if  $0<\gamma < 2\a$. Thus,   if $\gamma$ is chosen so that: $- \underline{\eps} + \gamma\sigma <0$, we get the result.
Eventually, for $\gamma \in (\nicefrac{\underline{\eps}}{\sigma '} , \nicefrac{\underline{\eps}}{\sigma })$, we obtain:
$$ u(x,t) \rightarrow 1 \mbox{ uniformly in } \{ \abs x \leq e^{\sigma t} \mbox{ as  } \a \rightarrow 1, t \rightarrow +\infty, t > \overline{C}\tau_{\a}\}.$$
\end{dem}

\section{Nondecreasing initial data}
The plan is similar to that of Section \ref{4}: first, we account for the linear propagation phase and then,  we describe the exponential propagation phase. Unfortunately, we cannot simply invoke Theorem $\ref{theo1}$ to prove $\ref{theo2}$: the computations are different, although related in spirit. 

Notice that the propagation exponent is (similarly to $\cite{JMRXC}$) strictly larger than in the compactly supported case.
\subsection{The linear propagation phase}
Recall that the problemn under consideration is:
\begin{equation} \label{systeme3}
\left\{
\begin{array}{rcl}
u_t +(-\partial_{xx})^{\a}u&=&u-u^2, \quad \R, t>0\\
u(x,0)&=&u_0(x), \hfill x \in \R\\
\end{array}
\right.
\end{equation}
where $u_0\in [0,1]$ is measurable, nondecreasing, such that $\underset{x \mapsto +\infty}{\lim}u_0(x)=1$ and for some constant $c$:
$$u_0(x) \leq c e^{-\abs x ^{\a}} \mbox{ for $x \in \R_-$ }.$$
\begin{lemme}\label{ite3}
For every $0<\sigma<2$ and $\a \in (\nicefrac{1}{2},1)$, there exist $\eps_0 \in (0,1)$ and $T_0 \geq 1$ depending only on $\sigma$ and $\eps_0$ for which the following holds.
Given $r_0 \leq -1$, $C$ independent of $\a$ and $\eps \in (0, \eps_0)$, let $\underline{u_0}=\eps \mathds{1}_{(r_0,+\infty)}$.
Then, the solution to \eqref{systeme} with initial condition $\underline{u_0}$ satisfies, for all $k\in \N$ such that $kT_0<\tau_{\a}$:
$$u(x,kT_0) \geq \eps \mbox{  for  } x \geq r_0 - k\sigma T_0^{\nicefrac{1}{\a}} .$$
\end{lemme}
\begin{dem}
For $k=0$, the result is obvious.\\
For $k=1$,Recall that for every $\delta \in (0,1)$, as long as $u\leq \delta$ we have :
$$ \underline{u}(x,t):=e^{(1-\delta)t} \int_{\R} \underline{u_0}(y) p(x-y,t) dy \leq u(x,t) \leq e^{t} \int_{\R} \underline{u_0}(y) p(x-y,t) dy=:\overline{u}(x,t).$$
Let $\delta >> \eps$, $T_0 > 0$ chosen so that:
\begin{equation}\label{T3}
\forall t\in (0,T_0), \forall x \in \R, \quad \overline{u}(x,t) \leq \delta.
\end{equation}
Note that the initial condition $u_0=\eps  \mathds{1}_{(r_0,+\infty)}$ is a nondecreasing function, so $u(\cdot,t)$ is nondecreasing for all $t>0$.
So it is sufficient to estimate: $\underset{ x \mapsto +\infty}{\lim}\overline{u}(x,t) \leq \delta$, to get \eqref{T3}. The inequalities on $p$ obtained in Proposition \ref{in1D} lead to, for $x$ large enough:
$$
\begin{disarray}{rcl}
\overline{u}(x,t) &\leq& e^t \int_{\R} \eps \mathds{1}_{(r_0,+\infty)}(y) p(x-y,t)dy  \\
&\leq& C\eps e^t \left(\int_{\abs{x-y} \leq \tild{C} t^{\nicefrac{1}{2\a}}} t^{-\nicefrac{1}{2\a}}dy+\right.\\
&& \hspace {3.5cm}\left.\int_{\abs{x-y} \geq \tild{C} t^{\nicefrac{1}{2\a}}} \frac{(1-\a)t^2}{\abs{x-y} ^{1+4\a}}+\frac{ \sin(\a\pi)t}{\abs{x-y} ^{1+2\a}} + \frac{e^{-\frac{\abs{x-y}}{4t}^{2\a}}}{\sqrt{t}\abs {x-y}^{1-\a}}dy\right) \\
&\leq&  C \eps e^t \left(1+(1-\a) +  \sin(\a\pi) + \int_{\abs z \geq 1} \frac{ e^{-\abs z^{2\a}}}{\abs z^{1-\a}} dz \right)\\
&\leq & B \eps e^t ,
\end{disarray}
$$
where $B$ is independent of $\a$. Thus, $T_0 = \ln\left( \dis{\frac{\delta}{B\eps}} \right)$ is smaller than $ \tau_{\a}$  (taking $1-\a$ smaller if necessary) and we have: 
$$\forall t\in (0,T_0), \forall x \in \R, \quad u(x,t) \leq \delta.$$
Note that we have the same time as in the case of a compactly supported initial datum.
Then, we look for $r_1< r_0<0$ for which we have: $u(x,T_0) \geq \eps,  \forall  x > r_1$. 
To find it, we look for $x_1<r_0$ such  that $\underline{u}$ is larger than $\eps$ for $ x$ larger than $x_1$.
First, let us notice that, for $y \in \R$, if $\abs{y-x} \geq \tild{C} T_0^{\nicefrac{1}{2\a}}$, where $\tild{C}^{2\a}\geq\dis{\frac{C (1-a)}{\sin(\a \pi)}}$, $\tild{C}$ independent of $\a$, then: 
\begin{equation} \label{maj} \frac{C(1-\a) T_0^2}{\abs{y-x}^{1+4\a}} \leq \frac{\sin(\a \pi) T_0}{ \abs{y-x}^{1+2\a}}.
\end{equation}
Next, we prove the existence of a constant $\overline{C} > 4^{\nicefrac{1}{2\a}}$ such that:
$$\tild{C}T_0^{\nicefrac{1}{2\a}} \leq  \abs{x_1-r_0} \leq \overline{C} T_0^{\nicefrac{1}{\a}}.$$
Indeed:
\begin{itemize}
\item if for all $\overline{C} >  4^{\nicefrac{1}{2\a}}$, $\abs{x_1-r_0} > \overline{C} T_0^{\nicefrac{1}{\a}}$, then $\overline{u}(x_1,T_0)$ is strictly smaller than $\eps$, which is impossible. Indeed:
$$
\begin{disarray}{rcl}
\overline{u}(x_1,T_0) &\leq& Ce^{T_0} \eps \left( \int_{ y > r_0} \frac{{T_0}^2  (1-\a)}{\abs{x_1-y}^{1+4\a}}+\frac{ T_0 \sin (\a \pi) }{\abs{x_1-y}^{1+2\a}} +\frac{e^{-\frac{\abs{y-x_1}}{4T_0}^{2\a}}}{\sqrt{T_0}\abs{x_1-y}^{1-\a}} dy \right)\\
&\leq&   \frac{C \eps {T_0}^2}{ \abs{x_1-r_0}^{4\a}} +   \frac{C\eps T_0}{ \abs{x_1-r_0}^{2\a}}  +C\eps  e^{T_0}  \int_{z \geq \frac{r_0-x_1}{(4T_0)^{\nicefrac{1}{2\a}}}} e^{- z^{2\a}}z^{\a-1}dz \\
&\leq&   \frac{C \eps }{T_0^2} +   \frac{C\eps}{T_0}  +\frac{C\eps  e^{T_0-\frac{\overline{C}^{2\a}}{4} T_0}}{T_0}  \\
&<& \eps, 
\end{disarray}
$$
 for $T_0$ large enough, since $e^{T_0} \leq C(1-\a)^{-1}$.
 So, there exists $\overline{C} >  4^{\nicefrac{1}{2\a}}$ such that $\abs{x_1-r_0} \leq \overline{C} T_0^{\nicefrac{1}{\a}}$.
\item if $ \abs{x_1-r_0} < \tild{C}T_0^{\nicefrac{1}{2\a}}$, then a sub solution is larger than $2 \eps$, which is impossible: indeed, since $u(\cdot,t)$ is nondecreasing for all $t>0$, we have, as in Lemma \ref{it}:
$$
\begin{disarray}{rcl}
u(x_1,T_0) &\geq& u(r_0-\tild{C}T_0^{\nicefrac{1}{2\a}},T_0)\\
&\geq&e^{(1-\delta)T_0} \int_{\R} \underline{u_0}(y) p(r_0-\tild{C}T_0^{\nicefrac{1}{2\a}}-y,T_0) dy\\
&\geq& e^{(1-\delta)T_0} \eps  p( \tild{C}T_0^{\nicefrac{1}{2\a}},T_0) \\
&\geq& 2 \eps, \\
\end{disarray}
$$
 for $T_0$ large enough.
\end{itemize}
Thus, we lay the emphasis on the Gaussian term, using \eqref{maj}, and for all $x$ such that:  $\tild{C} T_0^{\nicefrac{1}{2\a}} \leq \abs{x-r_0} \leq C T_0^{\nicefrac{1}{\a}} $, 
$$
\begin{disarray}{rcl}
u(x,T_0) &\geq& e^{(1-\delta)T_0} \int_{\R} p(x-y,T_0) \underline{u_0}(y) dy \\
&\geq&Ce^{(1-\delta)T_0}\eps\int_{\tiny{\begin{array}{l}\abs{y-x}\geq \tild{C} T_0^{\nicefrac{1}{2\a}} \\ y > r_0\end{array}}}\frac{e^{-\frac{\abs{y-x}}{4T_0}^{2\a}}}{\sqrt{ T_0}\abs{y-x}^{1-\a}} dy\\
&\geq&Ce^{(1-\delta)T_0}\eps\frac{e^{-\frac{\abs{r_0-x}}{4T_0}^{2\a}}}{\sqrt{T_0}\abs{r_0-x}^{1-\a}} \\
&\geq& C\eps e^{(1-\delta)T_0-\frac{\abs{r_0-x}}{4T_0}^{2\a}}T_0^{\nicefrac{1}{2}-\nicefrac{1}{\a}}.
\end{disarray}
$$
Let us define $x_1$ by: $$C\eps e^{(1-\delta)T_0-\frac{\abs{r_0-x_1}}{4T_0}^{2\a}}T_0^{\nicefrac{1}{2}-\nicefrac{1}{\a}}=\eps,$$
that is to say: $$ x_1= r_0-2^{\nicefrac{1}{\a}}T_0 ^{\nicefrac{1}{\a}}\left(1-\delta-\frac{1}{4T_0}\ln \left(CT_0^{\nicefrac{1}{\a}-\nicefrac{1}{2}}\right) \right)^{\nicefrac{1}{2\a}}.$$
Consequently, for $\sigma < 2$, we take $\delta$ small enough, $T_0$ large enough (but smaller than $\tau_{\a}$), and $\a$ close enough to $1$ so that:
$$\sigma < 2^{\nicefrac{1}{\a}} \left(1-\delta-\frac{1}{4T_0}\ln \left(CT_0^{\nicefrac{1}{\a}-\nicefrac{1}{2}}\right) \right)^{\nicefrac{1}{2\a}}.$$
Now, let us define $r_1:=r_0-\sigma T_0^{\nicefrac{1}{\a}}<r_0$ so that $ r_1>x_1 $. Thus since $u(\cdot,T_0)$ is a nondecreasing function:
$$ u(x,T_0) \geq u(x_1,T_0) \geq \eps, \ \forall  x \geq r_1.$$
Finally, $u(.,T_0) \geq \eps \mathds{1}_{(r_1,+\infty)}=: \underline{\underline{u_0}}$ and we can repeat the argument above, now with initial time $T_0$ and inital condition $\underline{\underline{u_0}}$ as long as $kT_0 < \tau_{\a}$ and get that:
$$u(x, kT_0)\geq \eps \mbox{ for } x \geq r_k,$$
for all $k \in \N$ satisfying $kT_0 < \tau_{\a}$, with: $r_k \geq r_0-\sigma k T_0^{\nicefrac{1}{\a}}.$
\end{dem}
\begin{cor} \label{coro3}
For every $0<\sigma < 2$ and $\a \in (\nicefrac{1}{2},1)$, let $T_0$ be as in Lemma \ref{ite3}. Then, for every $u_0\in [0,1]$  measurable, nondecreasing, such that $\underset{x \mapsto +\infty}{\lim}u_0(x)=1$, there exist $\eps \in (0,1)$ and $b>0$ such that 
$$u(x,t) \geq \eps, \mbox{ for }  T_0\leq t < \tau_{\a} \mbox{ and }  x \geq b-\sigma t^{\nicefrac{1}{\a}}.$$
\end{cor}
\begin{dem}
Let $\sigma \in (0,2).$
We have $u(\cdot,\tau)>0$ in $\R$ for all $\tau \in [\nicefrac{T_0}{2},\nicefrac{3T_0}{2}]$. Thus, there exists $\tild{\eps} \in (0,1)$ such that for all $\eps \in (0,\tild{\eps})$:
$$ u(\cdot,\tau)\geq {\eps} \mathds{1}_{(0, +\infty)}, \mbox{ in } \R$$
Let us define $\underline{u_0}$ by: $\underline{u_0}(y)={\eps} \mathds{1}_{(0, +\infty)}(y)$. Thus, $u(\cdot,\tau+t) \geq v(\cdot,t)$, for $t>0$, where $v$ is the solution to \eqref{systeme3} with initial condition $\underline{u_0}$ at time $\tau$. Next, we apply Lemma \ref{ite3} to $v$:  for all $k\in \N$ such that $kT_0 < \tau_{\a}$, we have:
$$v(x,kT_0) \geq \eps \mbox{  for  } x \geq  - k\sigma T_0^{\nicefrac{1}{\a}} .$$
Consequently:
For all $\tau \in  [\nicefrac{T_0}{2},\nicefrac{3T_0}{2}]$ , for all $k\in \N$ such that $\tau + kT_0< \tau_{\a}$:
$$u(x,\tau +kT_0) \geq \eps \mbox{  for  }  x \geq - k\sigma T_0^{\nicefrac{1}{\a}} .$$
Moreover, the set $\{\tau +kT_0, k \in \N , \tau \in  [\nicefrac{T_0}{2},\nicefrac{3T_0}{2}] \ | \  \tau+ kT_0 < \tau_{\a}\}$ covers all $[\nicefrac{T_0}{2}, \tau_{\a})$.
Let $ T_0 < t < \tau_{_a} $, then there exist $\tau  \in  [\nicefrac{T_0}{2},\nicefrac{3T_0}{2}]$ and $ k \in \N $ such that: $t=\tau + kT_0$.\\
Then for $\a$ close to $1$ enough, there exists a constant $b>0$  independent of $\a$ such that
 $$u(x,t)\geq \eps \mbox{ for }  x \geq b- \sigma t^{\nicefrac{1}{\a}}.$$
\end{dem}

Now, we can prove the first part of Theorem \ref{thm2}:
\begin{theo}\label{theo2}
Consider $u$ the solution to \eqref{systeme3}, with measurable and nondecreasing initial datum satisfying $u_0 \in [0,1]$, $\underset{x \mapsto +\infty}{\lim}u_0(x)=1$ and $$u_0(x) \leq c e^{-\abs x ^{\a}}, \mbox{ for x $\in \R_-$ },$$ for some constant $c$. Then:
\begin{itemize}
\item if $\sigma > 2$, then $u(x,t) \rightarrow 0$ uniformly in $\{  x \leq -\sigma t^{\nicefrac{1}{\a}}\}$ as $\a \rightarrow 1, t \rightarrow +\infty, t < \tau_{\a}$.
\item if $0< \sigma < 2$, then $u(x,t) \rightarrow 1$ uniformly in $\{  x \geq -\sigma t^{\nicefrac{1}{\a}}\}$ as $\a \rightarrow 1, t \rightarrow +\infty, t<\tau_{\a}$.
\end{itemize}
\end{theo}

\begin{dem}
We prove the first statement of the theorem.\\
Let $\sigma$ be such that $\sigma >2$, and $x \leq -\sigma t^{\nicefrac{1}{\a}}$. Recall that $u(x,t) \leq \overline{u}(x,t)$, where:
$$\overline{u}(x,t)=e^t \int_{\R} u_0(y) p(x-y,t) dy.$$
Let $\sigma' \in (2, \sigma)$. Using the fact: $e^{\tau_{\a}}=e^{\frac{\xi_{\a}^{2\a}}{4}}\underset{\a \rightarrow 1}{\sim}(1-\a)^{-1}$, for $t < \tau_{\a}$, we get:
$$
\begin{disarray}{rcl}
u(x,t) & \leq & Ce^t \left( \int_{\abs{x-y} \leq 1} \frac{u_0(y) }{t^{\nicefrac{1}{2\a}}}dy+\int_{\abs{x-y} \geq 1} \frac{t^2 (1-\a)u_0(y)}{\abs{x-y}^{1+4\a}}+ \frac{t\sin (\a \pi)u_0(y)}{\abs{x-y}^{1+2\a}}\right. \nonumber \\
&& \hspace{10cm}\left.+\frac{e^{-\frac{\abs{y-x}}{4t}^{2\a}} u_0(y)}{\sqrt{t}\abs{x-y}^{1-\a}} dy \right) \nonumber \\
&\leq &  C \left(  \int_{\abs{x-y} \leq 1}e^t e^{- \abs y ^{\a}}dy+\int_{y \leq x-1}e^{- \abs y ^{\a}}e^tdy  + \int_{x+1 \leq y \leq (\sigma '-\sigma) t^{\nicefrac{1}{\a}}}\frac{e^{- \abs y ^{\a}}t^2}{\abs{x-y}^{1+4\a}}dy  \right.\nonumber \\
&&   +\int_{x+1 \leq y \leq (\sigma '-\sigma) t^{\nicefrac{1}{\a}}}\frac{e^{- \abs y ^{\a}}t}{\abs{x-y}^{1+2\a}} +e^t\frac{e^{- \abs y ^{\a}}e^{-\frac{\abs{y-x}}{4t}^{2\a}} }{\abs{x-y}^{1-\a}} dy +\int_{ (\sigma '-\sigma) t^{\nicefrac{1}{\a}} \leq y }\frac{t^2}{\abs{x-y}^{1+4\a}}  dy \nonumber \\
&&\left.+ \int_{ (\sigma '-\sigma) t^{\nicefrac{1}{\a}} \leq y } \frac{t}{\abs{x-y}^{1+2\a}}+ e^t\frac{e^{- \abs y ^{\a}}e^{-\frac{\abs{y-x}}{4t}^{2\a}} }{\abs{x-y}^{1-\a}} dy \right)\\
&\leq &  C \left(  e^t e^{- (-1-x) ^{\a}} +e^t e^{- (1-x)^{\a}}+t^{\nicefrac{1}{\a}+1} e^{-(\sigma - \sigma ' )^{\a} t}+  t^{\nicefrac{1}{\a}}e^{-(\sigma - \sigma ' )^{\a} t} \nonumber\right. \\
&&   +\int_{x+1 \leq y \leq  (\sigma '-\sigma ) t^{\nicefrac{1}{\a}}}e^te^{- \frac{(-x )}{4t}^{2\a}}e^{\frac{((-x)^{\a} -2t)^2}{4t}}e^{-\frac{(y^{\a}-((-x)^{\a} -2t))}{4t}^2 } dy  \nonumber \\
&& + \int_{z \geq (\sigma '-\sigma) t^{\nicefrac{1}{\a}}-x \geq \sigma '  t^{\nicefrac{1}{\a}}}\frac{t^2}{\abs{z}^{1+4\a}}  +\frac{t}{\abs{z}^{1+2\a}} + e^{t-\frac{\abs{z}}{4t}^{2\a}}dz \nonumber \\
&\leq &  C \left(  e^{t-\frac {\sigma^{\a}}{2^{\a}}t}+t^{\nicefrac{1}{\a}+1} e^{-(\sigma - \sigma ' )^{\a} t}  +e^te^{t-(-x )^{\a}}\int_{\R}e^{- z ^2}dz  + \frac{1}{t^2}  +\frac{1}{t} + e^{t-\frac{\sigma '^{2\a}}{4} t} \right) \nonumber \\
&\leq &  C \left(  e^{t-\frac {\sigma^{\a}}{2^{\a}}t}+t^{\nicefrac{1}{\a}+1} e^{-(\sigma - \sigma ' )^{\a} t}  +e^{2t-\sigma^{\a}t} + e^{t-\frac{\sigma '^{2\a}}{4} t} \right) \nonumber \\
\end{disarray}
$$
Consequently, since $\sigma> \sigma '>2$, there exists $\a_2 \in (\a_1,1)$ such that $\sigma > 2^{\nicefrac{1}{\a_2}}$.
Thus: $(1-\frac{\sigma '^{2\a}}{4})<0, \mbox{ and } (2-\sigma^{\a})<0 , \ \forall \sigma \in (\a_2,1)$, and so:
$$ u(x,t) \rightarrow 0 \mbox{ uniformly in }\{ x \geq -\sigma t^{\nicefrac{1}{\a}}\} \mbox{ as }\a \rightarrow 1, t \rightarrow +\infty, t<\tau_{\a}.$$
Now we prove the second statement of the theorem:\\
Given $0<\sigma < 2$, take $\sigma' \in (\sigma,2)$, and apply Corollary \ref{coro3} with $\sigma$ replaced by $\sigma'$. Thus, we obtain:
$$-u \leq -\eps \mbox{ in } \omega:=\left\{ (x,t) \in \R \times \R^+ \mbox{ such that } T_0 \leq t< \tau_{\a} , \  x \geq b-\sigma' t^{\nicefrac{1}{\a}} \right\},$$
for some $b>0$.
Moreover: $(\partial_t + (-\Delta)^{\a})(1-u)=-u(1-u) \leq - \eps (1-u) \mbox{ in } \omega.$
Let $v$ be the solution to:
\begin{equation} 
\left\{
\begin{array}{rclc}
v_t +(-\Delta)^{\a}v&=&- \eps v,& \quad \R, t>0\\
v(y,T_0)&=&2+\frac{e^{\gamma \abs y^{\a}}}{D}  \mathds{1}_{\left\{ - c {\tau_{\a}}^{\nicefrac{1}{\a}}\leq y \leq 0\right\}}(y), &\: y \in \R,\\
\end{array}
\right.
\end{equation}
where $\gamma$ and $D$ are constants (depending on $\a$) chosen later,  $c$ is a constant  independent of $\a$ and strictly bigger than $\sigma'$.
There holds:
$$v(x,t)=e^{-\eps (t-T_0)} \left(2+ \dis{\int_{ - c {\tau_{\a}}^{\nicefrac{1}{\a}}\leq y \leq 0}\frac{ e^{\gamma \abs y^{\a}}}{D}p(x-y,t-T_0) dy} \right),$$
and $1-u\leq v  \mbox{ in } \omega$.

Now we want to apply Lemma 2.2 of \cite{JMRXC}, to have:
$$0\leq 1-u \leq v \mbox{ in } \R \times [T_0,\tau_{\a}).$$
Notice that we can apply the lemma since the assumption $\underset{x \mapsto +\infty}{\lim} u_0(x)=1$ leads to: $\underset{x \mapsto +\infty}{\lim} u(x,t)=1$, for all $t>0$.
Let $w:= 1-u-v$ with initial time $T_0$ and $ x \geq r(t):= b-\sigma' t^{\nicefrac{1}{\a}}$.
Let us verify the assumptions of the lemma:
\begin{itemize}
\item Initial condition: $w(.,T_0) \leq 0$ since $1-u \leq 1 \leq v \mbox{ for } t=T_0$
\item Condition outside $\omega$: let $\tau_{\a} > t \geq T_0$ and $ x \leq r(t)$. We have to verify that $w(x,t) \leq 0$, thus proving that $v(x,t) \geq 1$. We only have to consider the case $x \geq r(\tau_{\a})$. Indeed,  $\underset{x \mapsto -\infty}{\lim} v(x,T_0)=2$, so for all $t>T_0$ we have: $\underset{x \mapsto -\infty}{\lim} v(x,t)=2$. Consequently, for $\a$ close to $1$ enough, and $t\in [T_0, \tau_{\a})$:
$$ x \leq r(\tau_{\a})=b-\sigma' \tau_{\a}^{\nicefrac{1}{\a}} \Rightarrow v(x,t) \geq 1.$$
And for $\tau_{\a} > t \geq T_0$ and $ r(\tau_{\a}) \leq x \leq r(t)$:
$$
\begin{disarray}{rcl}
v(x,t) &\geq&C e^{-\eps (t-T_0)}\int_{\tiny{\begin{array}{l}\abs{y-x}\geq \tild{C}  (t-T_0)^{\nicefrac{1}{2\a}}  \\ - c {\tau_{\a}}^{\nicefrac{1}{\a}}\leq y \leq 0\end{array}}}\frac{e^{\gamma \abs y^{\a}}}{D} \frac{e^{-\frac{\abs{y-x}^{2\a}}{4 (t-T_0)}}}{\sqrt{ (t-T_0)}\abs{y-x}^{1-\a}} dy\\
&\geq& Ce^{-\eps (t-T_0)}\int_{\tiny{\begin{array}{l} \abs{z}\geq \nicefrac{\tild{C}}{2^{\nicefrac{1}{\a}}}  \\ \frac{x}{(4 (t-T_0))^{\nicefrac{1}{2\a}}}< z < \frac{x+ c {{\tau_{\a}}}^{\nicefrac{1}{\a}}}{(4 (t-T_0))^{\nicefrac{1}{2\a}}}\end{array}}}\frac{ e^{\gamma   \abs{ x-(4 (t-T_0))^{\nicefrac{1}{2\a}} z}^{\a}}}{D}\frac{e^{-\abs z^{2\a}}}{\abs{z}^{1-\a}} dz \\
&\geq& Ce^{-\eps (t-T_0)}\int_{\tiny{\begin{array}{l} \abs{z}\geq \nicefrac{\tild{C}}{2^{\nicefrac{1}{\a}}}  \\ \frac{x}{(4 (t-T_0))^{\nicefrac{1}{2\a}}}< z < \frac{b + (c-\sigma ') {{\tau_{\a}}}^{\nicefrac{1}{\a}}}{(4 (t-T_0))^{\nicefrac{1}{2\a}}}\end{array}}}\frac{ e^{\gamma   \abs{ x-(4 (t-T_0))^{\nicefrac{1}{2\a}} z}^{\a}}}{D}\frac{e^{-\abs z^{2\a}}}{\abs{z}^{1-\a}} dz\\
&\geq&Ce^{-\eps (t-T_0)} e^{\gamma (-x)^{\a}}
\end{disarray}
$$
We now choose $D$ so that the right hand side of the above inequality is larger than $ e^{(-\eps+\gamma    \sigma'^{\a}) (t-T_0)}$, and $\gamma$ such that: $-\eps + {\sigma'}^{\a} \gamma>0$.
\item Let $\tau_{\a} > t \geq T_0$ and $ x \geq r(t)$, then we have:
$$w_t(x,t) +(-\Delta)^{\a} w(x,t) \leq - \eps w(x,t),$$
and the last assumption is satisfied.
\end{itemize}
Thus: $w \leq 0$ in $\R \times [T_0, \tau_{\a})$, that is to say:
$$0\leq 1-u(x,t) \leq v(x,t)=e^{-\eps (t-T_0)} \left(2+ \dis{\int_{  - c {\tau_{\a}}^{\nicefrac{1}{\a}}\leq y \leq 0}\frac{ e^{\gamma \abs y^{\a}}}{D}p(x-y,t-T_0) dy} \right),$$
for all $ (x,t) \in \R \times [T_0, \tau_{\a}).$ Finally, we are going to prove that:
$v(x,t) \rightarrow 0$ uniformly in $\{  x \geq -\sigma t^{\nicefrac{1}{\a}}\}$ as $\a \rightarrow 1, t \rightarrow +\infty, t < \tau_{\a}$. For $t  < \tau_{\a}$ and $ x \geq -\sigma t^{\nicefrac{1}{\a}}$:\\
$$
\begin{disarray}{rcl}
v(x,t) &\leq &e^{-\eps(t-T_0)}\left(
2+\int_{\tiny{\begin{array}{l} - c {\tau_{\a}}^{\nicefrac{1}{\a}}\leq y \leq 0 \\ \abs{x-y} \leq [ \gamma (t-T_0)]^{\nicefrac{1}{\a}}\end{array}}}\frac{e^{\gamma \abs y^{\a}}}{D} p(x-y,t-T_0) dy\right.\\
&&\hspace{6cm}\left.+\int_{\tiny{\begin{array}{l} - c {\tau_{\a}}^{\nicefrac{1}{\a}}\leq y \leq 0 \\ \abs{x-y} \geq [  \gamma (t-T_0)]^{\nicefrac{1}{\a}}\end{array}}}\frac{e^{\gamma \abs y^{\a}}}{D} p(x-y, t-T_0) dy \right)\\
&\leq & Ce^{-\eps(t-T_0)}\left(2+\int_{\tiny{\begin{array}{l} - c {\tau_{\a}}^{\nicefrac{1}{\a}}\leq y \leq 0 \\ \abs{x-y} \leq [  \gamma (t-T_0)]^{\nicefrac{1}{\a}}\end{array}}}e^{\gamma \abs y^{\a}}dy \right. \\
&&+C e^{-\eps(t-T_0)} \left( \int_{\tiny{\begin{array}{l} - c {\tau_{\a}}^{\nicefrac{1}{\a}}\leq y \leq 0 \\ \abs{x-y} \geq [  \gamma (t-T_0)]^{\nicefrac{1}{\a}}\end{array}}}e^{\gamma \abs y^{\a}}\left(\frac{(t-T_0)^2  (1-\a)}{ \abs{x-y}^{1+4\a}}+   \frac{\sin(\a\pi) (t-T_0)}{\abs{x-y}^{1+2\a}}\right) \right.dy \\
&&\left.\hspace{6cm} +\int_{\tiny{\begin{array}{l} - c {\tau_{\a}}^{\nicefrac{1}{\a}}\leq y \leq 0\\ \abs{x-y} \geq [  \gamma (t-T_0)]^{\nicefrac{1}{\a}}\end{array}}}   \frac{e^{\gamma \abs y^{\a}-\frac{\abs{x-y}^{2\a}}{4 (t-T_0)}}}{\sqrt{ (t-T_0)}\abs{x-y}^{1-\a}}dy\right)\\
&\leq & Ce^{-\eps(t-T_0)}\left(2+ (t-T_0)^{\nicefrac{1}{\a}} e^{\gamma (-x) ^{\a}}e^{\gamma^2(t-T_0)} +  e^{\gamma c^{\a}{\tau_{\a}}}(1-\a)+ e^{\gamma c^{\a}{\tau_{\a}}} \sin(\a\pi) \right.\\
&& \left.+\int_{\tiny{\begin{array}{l} \abs z \geq[ \frac{\gamma}{2} ]^{\nicefrac{1}{\a}}(t-T_0)^{\nicefrac{1}{2\a}}\geq 1 \\ \frac{x}{(4 (t-T_0))^{\nicefrac{1}{2\a}}}< z < \frac{x+ c {{\tau_{\a}}}^{\nicefrac{1}{\a}}}{(4 (t-T_0))^{\nicefrac{1}{2\a}}}\end{array}}}  e^{\gamma  (-x)^{\a}}e^{2 \gamma \sqrt{(t-T_0)} z^{\a}-\abs{z}^{2\a}}dz\right)\\
&\leq& Ce^{-\eps(t-T_0)}+ C e^{-\eps(t-T_0)}   (t-T_0)^{\nicefrac{1}{\a}} e^{\gamma \sigma^{\a}t}e^{\gamma^2(t-T_0)} + Ce^{-\eps(t-T_0)} \left( e^{(-1+\gamma c^{\a}){\tau_{\a}}}\right.\\
&&\left. +e^{\gamma \sigma^{\a}t}\int_{\R}  e^{-( z^{\a} - \gamma \sqrt{(t-T_0)})^2}e^{\gamma^2(t-T_0)}dz \right) \\
&\leq& Ce^{-\eps(t-T_0)}+ Ce^{(-\eps+\sigma^{\a}\gamma +\gamma^2)(t-T_0)}   (t-T_0)^{\nicefrac{1}{2\a}}e^{\gamma \sigma^{\a}T_0} + Ce^{-\eps(t-T_0)} e^{(-1+\gamma c^{\a}){\tau_{\a}}} \\
&&+e^{(-\eps+\sigma^{\a}\gamma +\gamma^2)(t-T_0)} e^{\gamma \sigma^{\a}T_0}\\
\end{disarray}
$$
We have the result if $\gamma$ is chosen so that: $-\eps+\sigma^{\a}\gamma +\gamma^2<0 $, that is to say: $\gamma < \gamma_1= \frac{-\sigma^{\a}+\sqrt{\sigma^{2\a}+4 \eps}}{2} \underset{\eps \rightarrow 0}{\sim} \frac{\eps}{\sigma^{\a}}.$ With such a $\gamma$, we have for $\eps$ small enough: $-1+c^{\a}\gamma<0$.
Eventually, taking $\eps$ smaller if necessary, we have $\eps \sigma'^{-\a} < \gamma_1$, and consequently, we take: $ \eps \sigma'^{-\a} < \gamma < \gamma_1$ and we get:
$$u(x,t) \rightarrow 1 \mbox{ uniformly in } \{  x \geq - \sigma t^{\nicefrac{1}{\a}}\} \mbox{ as } \a \rightarrow 1, t \rightarrow +\infty, t<\tau_{\a}.$$
\end{dem}

\subsection{The exponential propagation phase}\label{propexp}
Let us now worry about the behaviour of \eqref{systeme3} for $t\geq \tau_{\a}$. The initial data is:
 $u(x, \tau_{\a})= \int_{\R}u_0(y) p(x-y,\tau_{\a})dy$.
The argument developped in Section \ref{4} holds and so the study of the level set $\{ x \in \R \ | \ u(x,t) \geq \underline{\eps} \}$ is based on the one of $\{ x \in \R \ | \ u(x,t) \geq \eps_{\a} \}$.
\begin{lemme}\label{ite4}
For every $0<\sigma<\frac{1}{2\a}$ and $\a \in (\nicefrac{1}{2},1)$, there exist $\eps_{0} \in (0,1)$, $T_{\a} \geq 1$ depending on $\sigma$ and $\eps_0$ and of order $\tau_{\a}$, and $\tild{\tau_{\a}}>0$ for which the following holds.
Given $r_0 < - 1$, $\underline{\eps} \in (0, \eps_{0})$, $\eps_{\a} \in (0, \underline{\eps})$, $a_{0,\a}$ be defined by $a_{0,\a} \abs {r_0}^{-2\a}=\eps_{\a}$, and let 
$$
\underline{u_{0,\a}}(x)= \left\{
\begin{array}{ll}
a_{0,\a}\abs{r_0}^{-2\a} & \mbox{ if $ x \leq r_0$} \\
a_{0,\a} \abs{x}^{-2\a} & \mbox{ if $x \geq r_0$} 
\end{array}
\right.
$$
Then, the solution to $u_t+(-\Delta)^{\a}u=u-u^2$ with initial condition $\underline{u_{0,\a}}$ and initial time $\tau_{\a}$ satisfies, for all $k\in \N$:
$$u(x, \tau_{\a}+\tild{\tau_{\a}}+kT_{\a} \geq \underline{ \eps} \  \mbox{  for  } \  x \geq (r_0 +M) e ^{\sigma k (\tild{\tau_{\a}}+T_0)},$$
where $M$ is defined in section \ref{inter}. (Recall that $M$ is large enough so that the principal Dirichlet eigenvalue of $(-\Delta)^{\a}-I$ in $B_M$ is negative.)
\end{lemme}
\begin{dem}
Let $\kappa > 0$ large enough, and $\a_1 \in (\nicefrac{1}{2},1)$ such that: $\eps_{\a} :=\sin(\a \pi)^{1+\kappa} < \eps, \ \forall \a \in (\a_1,1).$ 
Define $\delta_{\a}= \sqrt{\eps_{\a}}$ and $T_{\a} = \ln\left( \dis{\frac{\delta_{\a}}{B\eps_{\a}}}\right)$ so that (see \ref{T_0}):
\begin{equation}\label{T3}
\forall t\in (0,T_{\a}), \forall x \in \R^d, \quad u(x ,\tau_{\a}+t) \leq \delta_{\a}.
\end{equation}
Now, we prove the lemma:\\
\textbf{For $k=0$:} we have: $u(x,\tau_{\a}) \geq \eps_{\a} \mbox{ for } x \geq r_0.$\\
Then, $u$ is bigger than the solution $v$ to \eqref{v} with initial condition $\eps_{\a} \mathds{1}_{(r_0,r_0+M-1)}$ at time $\tau_{\a}$. 
Then, we use Proposition \ref{youpi} to the solution $v$, and we get the existence of $r \in (0,M)$ such that: $$u(r_0+r,\tau_{\a}+\tild{\tau_{\a}})\geq v(r_0+r,\tau_{\a}+\tild{\tau_{\a}})\geq \underline{\eps}.$$
Finally, since $r<M$ and $u(\cdot,t)$ is a nondecreasing function for all $t \geq 0$, we get the inequality:
$$u(x, \tau_{\a}+ \tild{\tau_{\a}}) \geq  \underline{\eps}, \ \forall  x \geq r_0+M.$$
\textbf{For $k=1$:} First, we look for $x_1< r_0$ such that: 
$$ u(x,\tau_{\a}+\tild{\tau_{\a}}+ T_0) \geq \sin(\a \pi)^{1+\kappa}=\eps_{\a}, \ \mbox{ for } x \geq x_1.$$ 
For every $\delta \in (0,1)$, $\delta >> \underline{\eps}$, we have, using \ref{T3} and for $x < r_0$:
$$
\begin{disarray}{rcl}
u(x,\tau_{\a}+T_{\a}) 
& \geq &Ce^{(1-\delta_{\a})T_{\a}}\int_{\tiny{\begin{array}{l} \abs{x-y}\geq \tild{C}T_{\a}^{\nicefrac{1}{2\a}}\\ y \leq x \end{array}}} \frac{T_{\a}\sin(\a\pi)\underline{u_{0,\a}}(y)}{ \abs{x-y}^{1+2\a}}dy \nonumber\\
& \geq &Ce^{(1-\delta_{\a}-\delta)T_{\a}}\frac{a_{0,\a}}{ \abs x^{2\a}} \int_{z \geq \tild{C}>1}\abs{z}^{-(1+2\a)}dz, \\
\end{disarray}
$$
since for $\kappa > \frac{2}{\delta}-1$, $e^{\delta T_{\a}}\sin(\a\pi) \geq 1$.
Let  us define $x_1<0$ by:
 $$ Ce^{(1-\delta_{\a}-\delta)T_{\a}}\frac{a_{0,\a}}{  \abs{x_1}^{2\a}} =\eps_{\a}$$
Since $a_{0,\a}= \eps_{\a} \abs{r_0}^{2\a}$, we get:
$$x_1=C r_0 e^{\frac{1-\delta_{\a}-\delta}{2\a}T_{\a}} .$$
Consequently, for each $\a < 1$, and $0 < \sigma < \frac{1}{2\a}$, we take$\delta$ and  $ \delta_{\a}$ such that:
\begin{equation}\label{delta2}
 \sigma <\frac{1-\delta_{\a}-\delta}{2\a}< \frac{1}{2\a}.
\end{equation}
Thus, taking $\a$ closer to $1$ if necessary, we have: $$Ce^{\frac{1-\delta_{\a}-\delta}{2\a}T_{\a}}  \geq e^{\sigma T_{\a}}$$
Now, let us define: $\overline{r_1}= r_0 e^{\sigma T_{\a}}$ so that $\overline{r_1} > x_1$, and  since $u(\cdot,t)$ is a nondecreasing function for all $t \geq 0$, we obtain:
$$u(x, \tau_{\a}+T_{\a}) \geq u(x_1, \tau_{\a}+T_{\a})= \eps_{\a }, \mbox{ for }x \geq \overline{ r_1}.$$
Furthermore, we have:
$u(x, \tau_{\a}+T_{\a}) \geq \dis{ \frac{a_{1,\a}}{\abs x^{2\a}}} \mbox{ for }  x \geq \overline{r_1},$ where  $ a_{1,\a} = \eps_{\a} \abs{x_1}^{2\a}$.
Thus:$$u(\cdot,\tau_{\a}+T_{\a}) \geq \underline{u_{1,\a}},$$
where $\underline{u_{1,\a}}$ is given by the same expression as $\underline{u_{0,\a}}$ with $(r_0,a_{0,\a})$ replaced by  $(\overline{r_1},a_{1,\a})$.
Finally, we use the case $k=0$ to make the connection with the level set $\{x \in \R^d \ | \ u(x,t)=\underline{\eps} \}$, replacing the initial time $\tau_{\a}$ by $\tau_{\a}+T_{\a}$ and the initial condition by $\underline{u_{1,\a}}$, and we use Proposition \ref{youpi} to get: 
$$u(x,\tau_{\a}+\tild{\tau_{\a}}+T_{\a}) \geq \underline{\eps}, \mbox{ for } x \geq r_1:=(r_0+M)e^{\sigma T_{\a}}.$$
We can repeat the argument above to get :
$$u(x, \tau_{\a}+\tild{\tau_{\a}}+kT_{\a}) \geq \eps_{\a}, \mbox{ for }x \geq r_k,$$
for all $k \in \N$, with $r_k \leq (r_0+M)e^{\sigma kT_{\a}}$.
\end{dem}
\begin{cor} \label{coro4}
For every $0<\sigma < \frac{1}{2\a}$ and $\a \in (\nicefrac{1}{2},1)$, let $T_{\a}$ and $\tild{\tau_{\a}}$ be as in Lemma \ref{ite4}. Then, for every $u_0\in [0,1]$  measurable, nondecreasing, such that $\underset{x \mapsto +\infty}{\lim}u_0(x)=1$ , there exist $\underline{\eps} \in (0,1)$, $\overline{C}>0$ a constant independent of $\a$, and $b_{\a}>0$ such that 
$$u(x,t) \geq \underline{\eps}, \mbox{ if } t \geq \overline{C} \tau_{\a} \mbox{ and } x \geq b_{\a} e^{\sigma t},$$
where $b_{\a}$ is proportional to $-e^{-\underline{C}\sigma \tau_{\a}}$, $ \underline{C}$ is a constant independent of $\a$ and strictly smaller than $\overline{C}$.
\end{cor}
\begin{dem}
By Lemma \ref{ite4}, there exists $\delta$ defined by \ref{delta2}.
Let $\kappa > \dis{\frac{2}{\delta}}+1$ and $\a_1 \in (\nicefrac{1}{2},1)$ such that: $\eps_{\a} :=\sin(\a \pi)^{1+\kappa} < \underline{ \eps}, \ \forall \a \in (\a_1,1).$ \\
Let $ \delta_{\a}$  to be chosen such that:  $\delta_{\a}=\sqrt{\eps_{\a}} $.
Using the proof done for $t<\tau_{\a}$, we get the existence of $\eps \in (0,1)$, for $\a$ closer to $1$ if necessary:
$$u(x, \tau_{\a}) \geq \eps \geq a_{0,\a} , \mbox{ for } x \geq  r_{\a},$$
where $-3\tau_{\a}^{\nicefrac{1}{\a}} < r_{\a} < -\tau_{\a}^{\nicefrac{1}{\a}}$, and $a_{0,\a}= \eps_{\a}  \abs{r_{\a}}^{2\a}$.
Thus $$u(\cdot, \tau_{\a}) \geq a_{0,\a} \mathds{1}_{(r_{\a},+\infty) } \mbox{ in } \R,$$
and $u(\cdot,  \tau_{\a}+t) \geq v(\cdot,t)$ for t>0, where $v$ is the solution to \eqref{systeme3} with initial condition $a_{0,\a} \mathds{1}_{(r_{\a},+\infty)}$. 
Recall that $\overline{T_0}$ the time before which the solution $u$ reaches $\delta$.
Moreover, we obtain, for $ x \leq r_{\a}\leq-\tild{C} t^{\nicefrac{1}{2\a}} $ for all $ t \in [\nicefrac{\overline{T_0}}{3},\nicefrac{\overline{T_0}}{3}+T_{\a}]$: (Recall that $T_{\a} < \frac{2}{3} \overline{T_0}$)

$$
\begin{disarray}{rcl}
u( x,  \tau_{\a}+t) &\geq & e^{(1-\delta)t}\int_{\R} p(x-y,t) a_{0,\a} \mathds{1}_{(r_{\a},+\infty)}(y) dy \\
&\geq&C \int_{\tiny{ \begin{array}{l} \abs{x-y} \geq \tild{C} t^{\nicefrac{1}{2\a}} \\  y \geq r_{\a} \end{array}}}e^{\frac{1-\delta}{3}\overline{T_0}} \frac{a_{0,\a} t \sin(\a \pi) } {\abs{x-y}^{1+2\a}}dy\\
&\geq&C\int_{\tiny{ \begin{array}{l} z \leq -\tild{C} t^{\nicefrac{1}{2\a}} \\ z \leq x-r_{\a} \end{array}}}\left(\frac{\delta}{B \eps_{\a}}\right)^{\frac{1-\delta}{3}}\frac{a_{0,\a} \sin(\a \pi)} { \abs{z}^{1+2\a}}dz\\
&\geq& C \frac{ a_{0,\a} } { \sin(\a \pi)^{(1+\kappa)(\frac{1-\delta}{3})-1}\abs{x}^{2\a}}\\
&\geq& C\frac{ a_{0,\a} } { \sin(\a \pi)^{\kappa_0}\abs{x}^{2\a}},
\end{disarray}
$$
where $\kappa_0  >0$.

Taking $\a$ closer to $1$ if necessary, we get, for $t \in [\nicefrac{\overline{T_0}}{3},\nicefrac{\overline{T_0}}{3}+T_{\a}]$:
$$u( x,  \tau_{\a}+t) \geq\dis{ \frac{a_{0,\a}}{\abs x^{1+2\a}}}, \mbox{ for }  x \leq r_{\a}.$$
As a consequence, using the fact $u(\cdot,t)$ is a non decreasing function for all $t \geq 0$:
 $$u( x,   \tau_{\a}+t) \geq u(r_{\a} , \tau_{\a} + t)\geq \eps_{\a}, \mbox{ for }  x \geq r_{\a}.$$
Finally:
$$ u(\cdot , \tau_{\a}+t) \geq \underline{ u_{0,\a}}, \ \ \forall  t \in [\nicefrac{\overline{T_0}}{3},\nicefrac{\overline{T_0}}{3}+T_{\a}],$$
where $\underline{u_{0,\a}}$ is the initial condition in Lemma \ref{ite4}.
Next, we can apply Lemma \ref{ite4} to the solution $ u(\cdot ,\cdot+\tau_0)$ for all  
$\tau_0 \in [  \nicefrac{\overline{T_0}}{3},  \nicefrac{\overline{T_0}}{3}+T_{a}]$.
Indeed, $\{ \tau_{\a}+\tau_0+\tild{\tau_{\a}}+kT_{\a}, \ k \in \N , \tau_0 \in  [  \nicefrac{\overline{T_0}}{3},  \nicefrac{\overline{T_0}}{3}+T_{\a} ] \}$ covers all $(\tild{\tau_{\a}}+\tau_{\a}+ \nicefrac{\overline{T_0}}{3}, + \infty)$.
Let $\overline{C} $ be a constant such that  $\tau_{\a}+ \nicefrac{\overline{T_0}}{3}+\tild {\tau_{\a}} \leq \overline{C}\tau_{\a}$. If $t \geq \overline{C}\tau_{\a}$, then there exist $\tau_0 \in  [ \nicefrac{\overline{T_0}}{3},\nicefrac{\overline{T_0}}{3}+T_{\a}]$ and $k \in \N$ such that: $$t= \tau_{\a}+\tau_0+kT_{\a}+\tild{\tau_{\a}}.$$
Then: $$u(x,t) \geq \underline{\eps}, \mbox{ if }  t \geq  \overline{C}\tau_{\a} \mbox{ and }  x \geq b_{\a} e^{\sigma t},$$
with : $b_{\a}=C \ e^{-\sigma \underline{C} \tau_{\a}} < 0$, where $C<0$ is independent of $\a$, and $\underline{C}$ is a constant independent of $\a$ such that $\tau_{\a}+ \nicefrac{\overline{T_0}}{3}+T_{\a} +\tild{\tau_{\a}}\geq \underline{C} \tau_{\a}$. Moreover we choose: $\underline{C} < \overline {C}$.
\end{dem}

Now, we can prove the second part of Theorem \ref{thm2}:
\begin{theo}
Under the assumptions of Theorem \ref{theo2}, there exists a constant $\overline{C} >0$ such that:
\begin{itemize}
\item if $\sigma > \dis{\frac{1}{2\a}}$, then $u(x,t) \rightarrow 0$ uniformly in $\{  x \leq -e^{\sigma t}\}$ as $\a \rightarrow 1, t \rightarrow +\infty, t > \overline{C} \tau_{\a}$.
\item if $0< \sigma < \dis{\frac{1}{2\a}}$, then $u(x,t) \rightarrow 1$ uniformly in $\{ x \geq -e^{\sigma t}\}$ as $\a \rightarrow 1, t \rightarrow +\infty, t >\overline{C}  \tau_{\a}$.
\end{itemize}
\end{theo}

Note that: $u_0(x) \leq c \abs x^{-1-2\a} \mbox{ for $ x \in \R_-$ and $\abs x$ large }$.\\
\begin{dem}
We prove the first statement of the theorem.
Let $\sigma$ be such that $\sigma >\dis{\frac{1}{2\a}}$, and $x$ such that $  x \leq -e^{\sigma t}$.
\begin{align*}
u(x,t)\leq& Ce^t \left( \int_{\abs{x-y} \leq 1} \frac{u_0(y) }{t^{\nicefrac{1}{2\a}}}dy+\int_{\abs{x-y} \geq 1} \frac{t^2 (1-\a)u_0(y)}{\abs{x-y}^{1+4\a}}+ \frac{t\sin (\a \pi)u_0(y)}{\abs{x-y}^{1+2\a}} +\frac{e^{-\frac{\abs{y-x}}{4t}^{2\a}} u_0(y)}{\sqrt{t}\abs{x-y}^{1-\a}} dy \right) \nonumber \\
\leq &  C \left( \int_{\abs{x-y} \leq 1} \frac{e^t}{\abs y^{1+2\a}t^{\nicefrac{1}{2\a}}}dy+\int_{y \leq x-1}\frac{e^t\left(t^2+t+\frac{1}{\sqrt{t}} \right)}{ \abs y^{1+2\a}} dy +\int_{x+1 \leq y \leq \nicefrac{x}{2} }\frac{e^t\left(t^2+t+\frac{1}{\sqrt{t}} \right)}{ \abs y^{1+2\a}} dy \right)  \nonumber \\
& + C\int_{y \geq \nicefrac{x}{2}}\frac{t^2 e^t }{\abs{x-y}^{1+4\a}}+ \frac{t e^t}{\abs{x-y}^{1+2\a}} +\frac{e^te^{-\frac{\abs{y-x}}{4t}^{2\a}} }{\sqrt{t}\abs{x-y}^{1-\a}} dy\\
\leq &  C \left( \frac{ e^t}{  (-1-x) ^{2\a}} +\frac {e^t\left(t^2+t+\frac{1}{\sqrt{t}} \right)}{(1-x)^{2\a}}  +\frac {e^t \left(t^2+t+\frac{1}{\sqrt{t}} \right) }{(-x)^{2\a}}+ \frac{t^2 e^t}{(-x)^{4\a}}+ \frac{te^t}{(-x)^{2\a}} +e^t e^{-\frac{(-x)}{4t}^{2\a}} \right)\\
\leq& C\left( e^{t-2 \a \sigma t}t^2+ t^2 e^{t-4 \a \sigma t}+e^{t- \frac{e^{\sigma^{\a}t^{\a}}}{2^\a}} \right)\\
\end{align*}
Hence, for $\sigma >\dis{\frac{1}{2\a}}$, we obtain:
$$ u(x,t) \rightarrow 0 \mbox{ uniformly in } \{  x \leq -e^{\sigma t}\} \mbox{ as  } \a \rightarrow 1, t \rightarrow +\infty, t > \tau_{\a}.$$
Now, we prove the second statement of the theorem. Using the proof done for $t < \tau_{\a}$, we know there exists $\eps \in (0,1)$ such that:$$u(\cdot, \tau_{\a}) \geq \eps \mathds{1}_{(r_{\a},+\infty)},$$
where $r_{\a}<0$ and bigger than $-2 \tau_{\a}^{\nicefrac{1}{\a}}$. 
As in the beginning of the proof of Corollary \ref{coro4}, we have, for all $t \in [ \nicefrac{\overline{T_0}}{3}, \nicefrac{ \overline{T_0}}{3}+T_{\a}]$:
$$u(x,t) \geq u_{0,\a} = \left\{ \begin{array}{rl}\frac{a_{0,\a}}{\abs x ^{2\a}}, &  x \leq r_{\a}\\
\eps_{\a}, &  x \geq r_{\a} \end{array} \right. . $$
Given $0<\sigma <\dis{\frac{1}{2\a}}$, let us take $\sigma' \in (\sigma,\dis{\frac{1}{2\a}})$ and  apply corollary \ref{coro4} with $\sigma$ replaced by $\sigma'$.
Thus, we obtain:
$$-u \leq -\underline{\eps} \mbox{ in } \omega :=\left\{ (x,t) \in  \R \times \R^+  \ | \ t\geq \overline{C}\tau_{\a},  x \geq b_{\a}e^{\sigma' t} \right\},$$
where  $b_{\a}=C \ e^{-\sigma' \underline{C} \tau_{\a}} < 0$.\\
Moreover: $(\partial_t + (-\Delta)^{\a})(1-u)=-u(1-u) \leq - \eps (1-u) \mbox{ in } \omega.$
Let $v$ be  the solution to:
\begin{equation} 
\left\{
\begin{array}{rclc}
v_t +(-\Delta)^{\a}v&=&- \underline{\eps} v,& \quad \R, t> \overline{C}\tau_{\a}\\
v(y, \overline{C}\tau_{\a})&=&1+\dis{\frac{e^{-\gamma\sigma\overline{C}\tau_{\a}}\abs y^{\gamma}}{D}}\mathds{1}_{\{y<0\}},& \: y \in \R\\
\end{array}
\right.
\end{equation}
where $\gamma \in (0,1)$  and $D$ are constants, independent of $\a$, chosen later.
This solution is given by: $$v(x,t)=e^{-\eps (t- \overline{C}\tau_{\a})} \left(1+ \dis{\int_{ y <0}\frac{e^{-\gamma\sigma\overline{C}\tau_{\a}} \abs y^{\gamma}}{D }p(x-y,t- \overline{C}\tau_{\a}) dy} \right),$$
and $1-u\leq v  \mbox{ in } \omega$.
To get:$$0\leq 1-u \leq v \mbox{ in } \R \times (\overline{C}\tau_{\a},+\infty).$$
Let us verify the assumptions of the Lemma 2.3 in \cite {JMRXC}.
Here, we use the fact $\underset{x \mapsto +\infty}{\lim} u_0 (x) = 1$ leads to $\underset{x \mapsto +\infty}{\lim} u (x,t) = 1, \forall t >0.$ \\
Let $w:= 1-u-v$ with initial time $\overline{C}\tau_{\a}$ and $ x \geq r(t):= b_{\a}e^{\sigma' t}$. Remember that in this case: $b_{\a}=(r_{\a}+M) \ e^{-\sigma' \underline{C} \tau_{\a}} < 0$.
\begin{itemize}
\item Initial datum: $w(.,\overline{C}\tau_{\a}) \leq 0$ since $1-u \leq 1 \leq v \mbox{ for } t=\overline{C}\tau_{\a}$
\item Condition outside $\omega$: let $ t \geq \overline{C}\tau_{\a}$ and $x \leq r(t)$ . We have to verify that $w(x,t) \leq 0$, proving that $v(x,t) \geq 1$.
 We use the same inequalities as before  taking $\overline{C}$ larger if necessary and using the fact that $\sigma < \sigma '$:
$$
\begin{disarray}{rcl}
v(x,t)&\geq& e^{-\underline{\eps}(t-\overline{C}\tau_{\a})}\int_{y <0}  \frac{ e^{-\gamma \sigma \overline{C}\tau_{\a}}\abs y^{\gamma}}{D} p(x-y, t-\overline{C}\tau_{\a})dy\\
&\geq& Ce^{-\underline{\eps}(t-\overline{C}\tau_{\a})}\int_{\tiny{\begin{array}{l}\abs{x-y}\geq \tild{C}(t-\overline{C}\tau_{\a})^{\nicefrac{1}{2\a}} \\ y<0 \end{array}}}\frac{\sin(\a\pi)(t-\overline{C}\tau_{\a})   e^{-\gamma \sigma \overline{C}\tau_{\a}} \abs y^{\gamma}}{D  \abs{x-y}^{1+2\a}} dy\\
&\geq& Ce^{-\underline{\eps}(t-\overline{C}\tau_{\a})}\int_{\tiny{\begin{array}{l}\abs{x+y}\geq \tild{C}(t-\overline{C}\tau_{\a})^{\nicefrac{1}{2\a}} \\ y> \abs x \end{array}}}\frac{(t-\overline{C}\tau_{\a})\sin(\a\pi)   e^{-\gamma \sigma \overline{C}\tau_{\a}}\abs y^{\gamma}}{D \abs{x+y}^{1+2\a}} dy\\
&\geq& Ce^{-\underline{\eps}(t-\overline{C}\tau_{\a})}  e^{-\gamma \sigma \overline{C}\tau_{\a}}\abs x^{\gamma}\sin(\a\pi)\int_{ \abs{z}\geq \tild{C}}\frac{ 1 }{D\abs{z}^{1+2\a}} dz\\
&\geq& \frac{C}{D \tild{C}^{2\a}}e^{-\underline{\eps}(t-\overline{C}\tau_{\a})}\abs{b_{\a}}^{\gamma} e^{-\gamma \sigma \overline{C}\tau_{\a}}e^{\sigma ' \gamma t}\sin(\a\pi)\\
&\geq& \frac{C}{D \tild{C}^{2\a}}  e^{(-\underline{\eps} +\gamma \sigma') ( t-\overline{C}\tau_{\a})} e^{\gamma \sigma ' \overline{C}\tau_{\a}-\gamma \sigma \overline{C}\tau_{\a}-\gamma  \sigma ' \underline{C}\tau_{\a}-\tau_{\a}}\\
&\geq& e^{(-\underline{\eps} +\gamma \sigma')(t-\overline{C}\tau_{\a})},
\end{disarray}
$$
where $D$ is chosen independent of $\a$ such that $D \tild{C}^{2\a}\leq C$.
Thus, if $\gamma$ satisfies: $-\underline{\eps} +\gamma \sigma' > 0$, then:
$$v(x,t) \geq 1 \geq 1-u(x,t) , \mbox{ for } x \leq r(t).$$
\item Let $ t \geq \overline{C}\tau_{\a}$ and $\abs x \leq r(t)$, then we have:
$$w_t(x,t) +(-\Delta)^{\a} w(x,t) \leq - \underline{\eps} w(x,t),$$
and the last assumption is satisfied.
\end{itemize}
So: $w \leq 0$ in $\R \times [\overline{C}\tau_{\a}, +\infty)$, that is to say:
$$0\leq 1-u(x,t) \leq v(x,t)=e^{-\eps (t-\overline{C}\tau_{\a})} \left(1+ \dis{\int_{  y <0 }}\frac{e^{-\gamma \sigma \overline{C}\tau_{\a}} \abs y^{\gamma}}{D}p(x-y,t-\overline{C}\tau_{\a}) dy \right),$$
for all $ (x,t) \in \R \times [\overline{C}\tau_{\a},+\infty).$
Finally, we are going to prove that: $v(x,t) \rightarrow 0$ uniformly in $\{  x \geq -e^{\sigma t}\}$ as $\a \rightarrow 1, t \rightarrow +\infty, t >\overline{C}\tau_{\a}$:
$$
\begin{disarray}{rcl}
v(x,t) &\leq& Ce^{-\underline{\eps} (t-\overline{C}\tau_{\a})} \left(1+\int_{\abs{x-y}\leq 1} \frac{e^{-\gamma \sigma \overline{C}\tau_{\a}} \abs y^{\gamma}}{D }dy\right.
+\int_{\tiny{\begin{array}{l} \abs{x-y} \geq 1 \\ y <0 \end{array}}} \frac{(1-\a)(t-\overline{C}\tau_{\a})^2 \abs y^{\gamma}}{De^{\gamma \sigma \overline{C}\tau_{\a}} \abs{x-y}^{1+4\a}}dy \\
&&+\int_{\tiny{\begin{array}{l} \abs{x-y} \geq 1 \\ y <0 \end{array}}} e^{-\gamma \sigma \overline{C}\tau_{\a}} \abs y^{\gamma} \frac{\sin(\a\pi) (t-\overline{C}\tau_{\a})}{D   \abs{x-y}^{1+2\a}} \left. +\frac{e^{-\frac{\abs{x-y}^{2\a}}{4(t-\overline{C}\tau_{\a})}} e^{-\gamma \sigma \overline{C}\tau_{\a}} \abs y^{\gamma}}{D \sqrt{ (t-\overline{C}\tau_{\a})} \abs{x-y}^{1-\a}}dy\right)\\
&\leq& C e^{-\underline{\eps} (t-\overline{C}\tau_{\a})} \left(1+ (- x+1)^{\gamma}+ \int_{\abs z \geq 1} \frac{(1-\a) (t-\overline{C}\tau_{\a})^2 e^{-\gamma \sigma \overline{C}\tau_{\a}} ( (- x)^{\gamma}+z^{\gamma})}{D \abs{z}^{1+4\a}} dz \right.\\
&&+\int_{\abs{z} \geq 1} \frac{\sin(\a \pi) (t-\overline{C}\tau_{\a})e^{-\gamma \sigma \overline{C}\tau_{\a}} ((-x)^{\gamma}+ z ^{\gamma})}{D  \abs{z}^{1+2\a}}
\left.+ \frac{e^{-\frac{\abs{z}^{2\a}}{4(t-\overline{C} \tau_{\a})}}e^{-\gamma \sigma \overline{C}\tau_{\a}} ((- x)^{\gamma}+  z ^{\gamma})}{D \abs{z}^{1-\a}}dz\right)\\
&\leq& Ce^{-\underline{\eps} (t-\overline{C}\tau_{\a})} \left(1+\int_{\abs z \geq 1} \frac{(t-\overline{C}\tau_{\a})^2z^{\gamma}}{ \abs z^{ 1 + 4\a} }dz\right.\left.+ \int_{\abs z \geq 1}\frac{(t-\overline{C}\tau_{\a})z^{\gamma}}{ \abs z^{ 1 + 2\a } }dz  +\int_{\R} z^{\gamma}e^{-\frac{\abs{z}^{2\a}}{4 (t-\overline{C}\tau_{\a})}}dz\right)\\
&&+Ce^{(-\underline{\eps}+\gamma \sigma)(t-\overline{C}\tau_{\a})} \left(1+\int_{\abs z \geq 1}  \abs z^{ - 1 - 4\a} dz+\int_{\abs z \geq 1}\abs z^{ - 1 - 2\a}dz  +\int_{\R} e^{-\frac{\abs{z}^{2\a}}{4 (t-\overline{C}\tau_{\a})}} dz\right).
\end{disarray}
$$
Notice that all the integrals converge if  $0<\gamma < 2\a$. Thus,   if $\gamma$ is chosen so that: $- \underline{\eps} + \gamma\sigma <0$, we get the result.
Eventually, for $\gamma \in(\nicefrac{\underline{\eps}}{\sigma '} , \nicefrac{\underline{\eps}}{\sigma })$, we obtain:
$$ u(x,t) \rightarrow 1 \mbox{ uniformly in } \{  x \geq -e^{\sigma t} \} \mbox{ as  } \a \rightarrow 1, t \rightarrow +\infty, t > \overline{C}\tau_{\a}.$$
\end{dem}

\bibliography{ref}

\begin{thebibliography}{10}

\bibitem{AW}
D.G. Aronson and H.F. Weinberger.
\newblock Multidimensional nonlinear diffusion arising in population genetics.
\newblock {\em Adv. Math. 30}, pages 33--76, 1978.

\bibitem{BHNI}
H.~Berestycki, F.~Hamel, and N.~Nadirashvili.
\newblock The speed of propagation for {KPP} type problems. {I} - {P}eriodic
  framework.
\newblock {\em J. Eur. Math. Soc.}, pages 173--213, 2005.

\bibitem{BHN}
H.~Berestycki, F.~Hamel, and N.~Nadirashvili.
\newblock The speed of propagation for {KPP} type problems. {II} - {G}eneral
  domains.
\newblock {\em J. Amer. Math. Soc. 23}, pages 1--34, 2010.

\bibitem{BRR}
H.~Berestycki, J.M. Roquejoffre, and L.~Rossi.
\newblock The periodic patch model for population dynamics with fractional
  diffusion.
\newblock {\em DCDS-S 4}, pages 1--13, 2011.

\bibitem{Blu}
R.~M. Blumenthal and R.~K. Getoor.
\newblock Some theorems on stable processes.
\newblock {\em Transactions of the American Mathematical Society 95}, pages
  263--273, 1960.

\bibitem{pre}
X.~Cabré and J.M. Roquejoffre.
\newblock Propagation de fronts dans les équations de {F}isher-{KPP} avec
  diffusion fractionnaire.
\newblock {\em Preprint}.

\bibitem{JMRXC}
X.~Cabré and J.M. Roquejoffre.
\newblock Front propagation in {F}isher-{KPP} equations with fractional
  diffusion.
\newblock {\em CRAS}, pages 1361--1366, 2009.

\bibitem{Cas}
D.~del Castillo-Negrete.
\newblock Truncation effects in superdiffusive front propagation with lévy
  flights.
\newblock {\em Physical Review 79}, pages 1--10, 2009.

\bibitem{Erd}
A.~Erd{\'e}lyi.
\newblock {\em Higher Transcendental Functions. {V}ol. {II}}.
\newblock McGraw-Hill Book Company, Inc., New York-Toronto-London, 1953.

\bibitem{Jimmy}
J.~Garnier.
\newblock Accelerating solutions in integro-differential equations.
\newblock {\em SIAM Journal on Mathematical Analysis}, 2010.

\bibitem{HR}
F.~Hamel and L.~Roques.
\newblock Fast propagation for {KPP} equations with slowly decaying initial
  conditions.
\newblock {\em J. Differential Equations}, pages 1726--1745, 2010.

\bibitem{Kolmo}
A.N. Kolmogorov, I.G. Petrovskii, and N.S. Piskunov.
\newblock Etude de l'équation de diffusion avec accroissement de la quantité de
  matière, et son application à un problème biologique.
\newblock {\em Bjul. Moskowskogo Gos. Univ. 17}, pages 1--26, 1937.

\bibitem{Kolo}
V.~Kolokoltsov.
\newblock Symmetric stable laws asn stable-like jump-diffusions.
\newblock {\em London Math. Soc. 80}, pages 725--768, 2000.

\bibitem{vulpi}
R.~Mancinelli, D.~Vergni, and A.~Vulpiani.
\newblock Front propagation in reactive systems with anomalous diffusion.
\newblock {\em Phys. D 185}, pages 175--195, 2003.

\bibitem{Polya}
G.~Polya.
\newblock On the zeros of an integral function represented by {F}ourier's
  integral.
\newblock {\em Messenger of Math. 52}, pages 185 --188, 1923.

\bibitem{sch}
L.~Schwartz.
\newblock {\em M\'ethodes math\'ematiques pour les sciences physiques}.
\newblock Enseignement des Sciences. Hermann, Paris, 1961.

\bibitem{uchi}
K.~Uchiyama.
\newblock The behavior of solutions of the equation of
  {K}olmogorov-{P}etrovsky-{P}iskunov.
\newblock {\em Proc. Japan Acad. Ser. A Math. Sci.}, pages 225--228, 1977.

\bibitem{Wein}
H.F. Weinberger.
\newblock On spreading speeds and traveling waves for growth and migration
  models in a periodic habitat.
\newblock {\em J. Math. Biol}, pages 511--548, 2002.

\end{thebibliography}

\end{document}